\documentclass[reqno]{amsart}
\newtheorem{thm}{Theorem}[section]

\newtheorem{prop}[thm]{Proposition}
\newtheorem{lem}[thm]{Lemma}
\newtheorem{cor}[thm]{Corollary}
\newtheorem{exm}[thm]{Example}
\usepackage{amsmath,amssymb,amscd,epsfig}
\usepackage{color, graphics}
\numberwithin{equation}{section}

\newcommand{\ld}{\lambda}

\newcommand{\Z}{\mathbb Z}

\begin{document}
\title[Generalized Burnside-Grothendieck ring functor]
{Generalized Burnside-Grothendieck ring functor
and aperiodic ring functor associated with profinite groups}
\author{YOUNG-TAK OH}
\address{Korea Institute for Advanced Study\\
         207-43 Cheongryangri-dong, Dongdaemun-gu\\
         Seoul 130-722, Korea}
\email{ohyt@kias.re.kr} \maketitle \baselineskip=12pt

\begin{abstract}
For every profinite group $G$, we construct two covariant functors
$\Delta_G$ and ${\bf {\mathcal {AP}}}_G$ from the category of
commutative rings with identity to itself, and show that indeed
they are equivalent to the functor $W_G$ introduced in [A. Dress
and C. Siebeneicher, The Burnside ring of profinite groups and the
Witt vectors construction, {\it Adv. Math.} {\bf{70}} (1988),
87-132]. We call $\Delta_G$ the generalized Burnside-Grothendieck
ring functor and ${\bf {\mathcal {AP}}}_G$ the aperiodic ring
functor (associated with $G$). In case $G$ is abelian, we also
construct another functor ${\bf Ap}_G$ from the category of
commutative rings with identity to itself as a generalization of
the functor ${\bf Ap}$ introduced in [K. Varadarajan, K. Wehrhahn,
Aperiodic rings, necklace rings, and Witt vectors, {\it Adv.
Math.} {\bf 81} (1990), 1-29]. Finally it is shown that there
exist $q$-analogues of these functors (i.e, $W_G,\,\Delta_G,\,
{\bf {\mathcal {AP}}}_G$, and ${\bf Ap}_G$) in case $G=\hat C$ the
profinite completion of the multiplicative infinite cyclic group.

\end{abstract}

\section{Introduction}
\renewcommand{\thefootnote}{}
\footnotetext[1]{%
\renewcommand{\baselinestretch}{1.2}\selectfont
This research was supported by KOSEF Grant \# R01-2003-000-10012-0.
\hfill \break
MSC : 11F03,11F22,17B70.
\hfill \break
Keywords: Necklace ring, Witt-Burnside ring, Burnside-Grothendieck ring}

Since E. Witt introduced it in \cite{WI},
the universal ring of Witt vectors
has attracted many mathematicians until now
due to its interesting and peculiar ring structure.
The universal ring of Witt vectors $W(A)$
associated with a commutative ring $A$ with $1\ne 0$
can be characterized by the following properties:

(W1) As a set, it is $A^{\mathbb N}.$

(W2) For any ring homomorphism  $f:A\to B$, the map
$W(f):{\bf a}\mapsto (f(a_n))_{n\ge 1}$ is a ring homomorphism
for ${\bf a}=(a_n)_{n\ge 1}$.

(W3) The maps $w_n: W(A)\to A$ defined by
$${\bf a} \mapsto \sum_{d|n}da_d^{\frac dn}
\text{  for }{\bf a}=(a_n)_{n\ge 1}$$
are ring homomorphisms.

In \cite {MR} Metropolis and Rota gave a new interpretation
of the ring of Witt vectors using the concept of so called {\it `string'},
and using it they found an algebraic structure isomorphic to $W(A)$
for integral domains $A$ of characteristic zero.
They called this algebraic structure {\it necklace ring
${\bf Nr}(A)$ over $A.$}
The underlying set of ${\bf Nr}(A)$ is $A^{\mathbb N}.$
Addition is defined componentwise and multiplication is defined by
\begin{equation*}
({\bf bc})_n= \sum_{[i,j]=n}(i,j)b_ic_j
\end{equation*}
for ${\bf b}=(b_n)_{n\ge 1},\,{\bf c}=(c_n)_{n\ge 1} \in {\bf Nr}(A).$
In particular, in \cite{DS}, the necklace ring ${\bf Nr}(\mathbb Z)$
over $\mathbb Z$ was interpreted as the Burnside-Grothendieck ring
$\hat \Omega(C)$ of isomorphism classes
of almost finite cyclic sets.
Here $C$ denotes the multiplicative infinite cyclic group.

On the other hand, in \cite{DS2} Dress and Siebeneicher showed
that the classical construction of Witt vectors can be explained
as a special case of the construction of a certain functor using
group-theoretical language. In detail, they proved that there
exists a covariant functor $W_G$ from the category of commutative
rings into itself associated with an arbitrary profinite group
$G$, and that the classical construction of Witt vectors can be
recovered by considering the case $G=\hat C$, the profinite
completion of $C.$ Furthermore they also proved that if $A=\mathbb
Z$ the corresponding Witt-Burnside ring $W_G(\mathbb Z)$ is
isomorphic to the Burnside-Grothendieck ring $\hat \Omega(G)$ of
isomorphism classes of almost finite $G$-spaces. In view of Dress
and Siebeneicher's work it would be quite natural to question
whether there exists a functor equivalent to $W_G$ which coincides
with $\hat \Omega(G)$ when $A=\mathbb Z$. In \cite{O,OH} the
author gave an affirmative answer for a certain class of
commutative rings (more correctly, the category of {\it special
$\ld$-rings}), and in \cite{Br} M. Brun showed that if $G$ is a
finite group then the functor $W_G$ is equivalent to the left
adjoint of the algebraic functor from the category of $G$-Tambara
functors to the category of commutative rings with an action of
$G$. In this paper we remove such a restriction and construct two
functors $\Delta_G$ and ${\bf {\mathcal {AP}}}_G$ from the
category of commutative rings with identity to itself in a purely
algebraic manner. Actually these two functors will turn out
to be equivalent to $W_G$.

Surprisingly these functors have
$q$-analogues in some cases. More precisely, in \cite{L} C. Lenart
showed that there exist well-defined $q$-deformations of
$W_G(R)$, ${\bf Nr}_G(R)$, and ${\bf {\mathcal {AP}}}_G(R)$ when
$G=\hat C$ for any integer $q$ and for any commutative
torsion-free ring $R$ with identity.
We observe that indeed his results work for every commutative ring $R$
with identity.
Based on this observation we construct $q$-analogues $W_G^q(R)$, ${\bf
Nr}_G^q(R)$, $\Delta_G^q(R)$, and ${\bf {\mathcal {AP}}}_G^q(R)$
when $G=\hat C$ for every commutative ring $R$ with identity, and
also show that their construction is functorial.

This paper is organized as follows.
In Section 2, we introduce the basic definitions
and notation needed throughout this paper.
In Section 3 and 4 we construct the functors $\Delta_G$,
${\bf {\mathcal {AP}}}_G$
and show that they are equivalent to $W_G$.
One of many important properties of
$W_G(R)$, $\Delta_G(R)$, and ${\bf {\mathcal {AP}}}_G(R)$
is that they always come equipped with two families of additive homomorphisms.
One is inductions and the other is restrictions.
In fact inductions are additive,
whereas restrictions are ring homomorphisms.
In Section 5 it is shown that
natural equivalences among $W_G$, $\Delta_G$, and ${\bf {\mathcal {AP}}}_G$
are compatible with inductions and restrictions.
Finally in Section 6 we deal with
$q$-analogues of these functors in case $G=\hat C.$
\vskip 3mm

\section{Preliminaries}
Throughout this paper the rings $R$ we consider
will be commutative associative rings with $1_R\ne 0$,
and the subrings will have $1_R$.
We begin with the basic definitions and notation
on the covariant functor $W_G$ introduced by Dress
and Siebeneicher.
For complete information we refer to \cite{DS2}.
Let $G$ be an arbitrary profinite group.
For any $G$-space $X$ and any subgroup $U$
of $G$ define $\varphi_U(X)$ to be the cardinality
of the set $X^U$ of $U$-invariant
elements of $X$ and let $G/U$ denote the $G$-space of
left cosets of $U$ in $G.$
With this notation Dress and Siebeneicher showed that
there exists a unique covariant functor $ W_G$
from the category of commutative rings into itself
satisfying the following three conditions:

(${W_G}$1) For any commutative ring $A$ the ring
$W_G(A)$ coincides, as a set, with the set
$A^{\underline {\mathcal O(G)}}$ of all maps from the set $\mathcal O(G)$
of open subgroups of $G$ into the ring $A$ which are constant on
conjugacy classes.

(${W_G}$2) For every ring homomorphism $h:A\to B$
and every $\alpha \in  W_G(A)$ one has
$ W_G(h)(\alpha)=h\circ \alpha.$

(${W_G}$3) For any open subgroup $U$ of $G$ the family of maps
$$\phi_U^A:  W_G(A) \to A$$
defined by
\begin{equation*}
\alpha \mapsto \underset{U\lesssim V\leqslant G}{{\sum}'}
\varphi_U(G/V)\cdot\alpha(V)^{(V:U)}
\end{equation*}
provides a natural transformation from the functor $ W_G$ into the identity.
Here $U\lesssim V$ means that the open subgroup $U$ of $G$ is subconjugate
to $V,$ i.e., there exists some $g\in G$ with $U\leqslant gVg^{-1},\,(V:U)$
means the index of $U$ in $gVg^{-1}$ and the symbol ``$\sum'$"
(and ``${\prod}'$", too)
is meant to indicate that for each conjugacy class of subgroups $V$ with
$U\lesssim V$ exactly one summand has to be taken.

Much effort has been made to find a functor equivalent to $W_G$
(for example, see \cite{Br,Gr,O,OH,VW}),
and in \cite{O,OH},
the author succeeded in getting such a functor ${\bf Nr}_G$
when we view $W_G$ as a functor from the category of special $\ld$-rings
to the category of commutative rings.
Let us recall its definition briefly.
For arbitrary special $\ld$-rings $R$,
{\it the necklace ring  ${\bf Nr}_G(R)$ of $G$ over $R$}
is defined as follows:

(${\bf Nr}_G 1)$ Its underlying set is
$$\underset{U\leqslant G,\, \text{open}}{{\prod}'}R.$$

(${\bf Nr}_G 2)$ Addition is defined componentwise.

(${\bf Nr}_G 3)$ Multiplication is defined so that
the $U$-th component of the product of two
sequences ${\bf x}=(x_V)'_{V\leqslant G}$ and
${\bf y}=(y_W)'_{W \leqslant G}$ is given by
\begin{equation*}
({\bf x\cdot y})_U=\left(\underset{V,W}{{\sum}'}\sum_{VgW\subseteq G \atop Z(g,V,W)}
\Psi^{(V:Z(g,V,W))}(b_V)\Psi^{(W:Z(g,V,W))}(c_W)\right),
\end{equation*}
where the sum is over
$$Z(g,V,W):=V\cap gWg^{-1}$$ which are conjugate to $U$ in $G$,
and $\Psi^n,\,n\ge 1$ represents the $n$-th Adams operation.

In proving that $W_G(R)$ is isomorphic to ${\bf Nr}_G(R)$
two families of maps play a crucial role.
One is {\it exponential} maps $\tau^U(R)$ and
the other is inductions $\text{Ind}_U^G(R),$
where $U$ is any open subgroup of $G$.

\vskip 2mm
\noindent{\bf Caution.}
In \cite {O,OH}
the notation $\tau^U_R$ has been used instead of $\tau^U(R).$
However, in this paper, the notation $\tau^U_R$ will denote the map
$\tau_R$ in \cite {O,OH}.

\vskip 2mm
In detail, the map
$\tau^U(R):R\to  {\bf Nr}_U(R)$ is defined by
\begin{equation*}
r\mapsto
\left({\bf M}_U(r,V)\right)'_{V\leqslant U}.
\end{equation*}
Recall that for any profinite group $G$ exponential maps are given by
\begin{equation*}
\tau^G(R)(r)=\underset{V\leqslant G}{{\sum}'}{\bf M}_G(r,V)\,[G/V].
\end{equation*}
Here the coefficients ${\bf M}_G(r,V)$ are determined in the following way.
First we write $r$ as a sum of one-dimensional elements, say,
$r_1+r_2+\cdots+r_m.$
Now, consider the set of continuous maps from $G$ to the topological space
${\mathfrak r}:=\{r_1,r_2,\cdots,r_m\}$
with regard to the discrete topology with trivial $G$-action.
It is well-known that this set
becomes a $G$-space
with regard to the compact-open topology
via the following standard $G$-action
$$(g\cdot f)(x)=f(g^{-1}\cdot x).$$
Consider a disjoint union of $G$-orbits of this set, say,
$\underset{h}{\dot\bigcup}
\,G\cdot h,$
where $h$ runs through a system of representatives of this decomposition.
After writing
$G/G_h
=\underset{1\le i \le (G:G_h)}{\dot\bigcup}w_iG_h,$
where $G_h$ represents the isotropy subgroup of $h$,
we let
$$[h]:=\prod_{i=1}^{ (G:G_h)}h(w_i).$$
Clearly this is well-defined since $h$ is $G_h$-invariant.
With this notation, ${\bf M}_G(r,V)$ is given by
$\sum_h [h],$
where $h$ is taken over the representatives such that
$Gh$ is isomorphic to $G/V$.
Note that the $Gh$ is isomorphic to $G/V$ if and only if
$G_h$ is conjugate to $V$.
Indeed
$G$ acts on $h$ freely modulo $G_h$, that is,
$w_i\cdot h=w_j\cdot h \Rightarrow i=j.$
For example, if $G=\hat C$,
$${\bf M}(r,n):={\bf M}_{\hat C}(r,\hat C^n)
=\dfrac 1n \sum_{d|n}\mu(d)\Psi^d(r)^{\frac nd},\,\,n\ge 1$$
(see \cite{OH}).
The notation $\mu$ represents the classical M\"obius inversion function.
On the other hand, inductions
$\text{Ind}_U^G(R):{\bf Nr}_U(R)\to {\bf Nr}_G(R)$ are defined by
\begin{equation*}
(x_V)'_{V\leqslant U} \mapsto  (y_W)'_{W\leqslant G},
\end{equation*}
where $y_W$ is the sum of $x_V$'s such that $V$ is conjugate to $W$ in $G$.
It is clear that $\text{\rm Ind}_U^G(R)$ is additive.

\begin{prop} {\rm (Oh \cite{O,OH})}\label{multiplicity of teich map}
For every open subgroup $U$ of $G$ and every special $\ld$-ring $R$,
$\tau^U(R)$ is multiplicative.
\end{prop}

Using inductions $\text{Ind}_U^G(R)$ and
exponential maps $\tau^U(R)$ for all open subgroups of $G$
simultaneously, the author constructed a map, called
{\it $R$-Teichm\"{u}ller map}, as follows:
\begin{align*}
&\tau_R^G:W_G(R)\to {\bf Nr}_G(R)\\
&\alpha\mapsto {\underset{U}{\sum}'}\,\text{Ind}_U^G(R)
\circ \tau^U(R)(\alpha(U)).
\end{align*}

\begin{prop} {\rm (Oh, Theorem 3.3 \cite{OH})} \label{ringhomo of teichmuller}
For every special $\ld$-ring $R$, $\tau_R^G$ is a ring isomorphism.
Moreover the following diagram
\begin{equation}
\begin{picture}(360,70)
\put(100,60){$W_G(R)$}
\put(135,63){\vector(1,0){70}}
\put(210,60){${\bf Nr}_G(R)$}
\put(155,0){$R^{\underline {\mathcal O(G)}}$}
\put(215,55){\vector(-1,-1){43}}
\put(120,55){\vector(1,-1) {43}}
\put(165,67){$\tau_R^G$}
\put(204,35){$\tilde\varphi_R$}
\put(122,34){$\Phi_R$}
\put(163,38){$\curvearrowright$}
\end{picture}
\end{equation}
is commutative.
Here
$$\Phi_R=\underset{U\leqslant G,\, \text{\rm open}}{{\prod}'} \phi_U^R$$
and
$$\tilde \varphi_R
=\underset{U\leqslant G,\, \text{\rm open}}{{\prod}'} \tilde\varphi_U^R\,,$$
where the map $\tilde\varphi_U^R:{\bf Nr}_G(R)
\to R$ is given by
\begin{equation*}
(x_V)'_V\mapsto
\underset{ U\lesssim V \leqslant G}{{\sum}'}
\varphi_U(G/V)\Psi^{(V:U)}(x_V)\,.
\end{equation*}
\end{prop}

We close this section by providing examples
associated with two interesting special $\ld$-rings.
\begin {exm}{\rm
(a) Letting $R:=\mathbb C[a_1,a_2,\cdots],$
$R$ has a special $\ld$-ring structure associated with
the Adams operators
$\Psi^n(a_m):=a_{nm},\, m,n \ge 1.$
Let $G$ be the profinite completion $\hat C$ of the multiplicative
infinite cyclic group $C$.
Then, for ${\bf a}=(a_1,a_1^2,a_1^3,\cdots)$, we have
\begin{equation}\label{commuting diagram*}
\begin{picture}(360,70)
\put(100,60){$(a_1,0,0,\cdots)$}
\put(155,63){\vector(1,0){30}}
\put(155,61){\line(0,1){4}}
\put(190,60)
{$\left(\dfrac 1n \sum_{d|n}\mu(d)a_d^{\frac nd}\right)_{n\in \mathbb N}$}
\put(150,0){$(a_1,a_1^2,a_1^3,\cdots)$}
\put(225,53){\vector(-1,-1){43}}
\put(223,54){\line(1,0) {4}}
\put(120,53){\vector(1,-1) {43}}
\put(118,54){\line(1,0) {4}}
\put(163,67){$\tau_R^{\hat C}$}
\put(215,35){$\tilde \varphi_R$}
\put(121,34){$\Phi_R$}
\put(170,35){$\curvearrowright$}
\end{picture}
\end{equation}
Here the term
$$\dfrac 1n \sum_{d|n}\mu(d)a_d^{\frac nd}$$
can be found as
the cyclic index $Z_{C_l}(\alpha_l)$, where
$C_l$ is the subgroup of the symmetric group $S_l$ generated by $(1,2,\cdots,l)$ and $\alpha_l$
is the linear character defined by
$\alpha_l(1,2,\cdots,l)=e^{2\pi i/l}$ (see \cite{JK}).

(b) Letting $R=\mathbb C[[z]]$,
$R$ has a special $\ld$-ring structure associated with
the Adams operators
$\Psi^n(z):=z^n,\,n \ge 0.$
Then, for a formal power series $f(z)\in R$, we have
$${\bf M}(f(z),n)=\dfrac 1n \sum_{d|n}\mu(d)f(z^d)^{\frac nd}.$$
The coefficients $m_f(i,n)$ defined by
${\bf M}(f(z),n)=:\sum_{i=0}^{\infty}m_f(i,n)z^i$
have many interesting properties (see \cite {MO}),
and many of which can be easily proved under our $\ld$-ring notation
(see [Section 3 \cite{OH}]).}
\end{exm}

\section{ Generalized Burnside-Grothendieck ring functor $\Delta_G$}
In this section, given any profinite group $G$,
we shall introduce a covariant functor $\Delta_G$
from the category of commutative rings to itself.
Indeed, it will turn out to be equivalent with
the functor $W_G$ in \cite{DS2}.
To do this we introduce a functor ${\bf NR}_G$
from the category of commutative rings to itself.
Let $x_U$ and $y_U$ be indeterminates where $U$
ranges over every open subgroup of $G$.
We consider the following system of equations:
for every open subgroup $U$ of $G$
it holds that
\begin{equation*}
\underset{ U\lesssim V \leqslant G}{{\sum}'}
\varphi_U(G/V)\,s_V
=\underset{ U\lesssim V \leqslant G}{{\sum}'}
\varphi_U(G/V)\,(x_V+y_V),
\end{equation*}
and
\begin{equation*}
\begin{aligned}
\underset{ U\lesssim V \leqslant G}{{\sum}'}
\varphi_U(G/V)\,p_V
=\underset{ U\lesssim V \leqslant G}{{\sum}'}
\varphi_U(G/V)\,\,x_V\cdot
\underset{ U\lesssim V \leqslant G}{{\sum}'}
\varphi_U(G/V)\,\,y_V.
\end{aligned}
\end{equation*}
It is clear that $s_G=x_G+y_G$ and $p_G=x_G\cdot y_G$.
Solving $s_U$ and $p_U$ inductively,
one can show that for every open subgroup $U$ of $G$,
\begin{equation}\label{operation of necklace--equation 1}
\begin{cases}
s_U=x_U+y_U\\
p_U=\underset{V,W}{{\sum}'}
\underset{g}{{\sum}}
\,\,x_V \cdot y_W,
\end{cases}
\end{equation}
where $g$ runs through a system of representations
of the double cosets $VgW \subseteq G$ such that
$Z(g,V,W)$ is conjugate to $U$ in $G$.
Indeed this observation enables us to state the following theorem.

\begin{thm}\label{def of the functor:NECKLACE RING}
Associated with every profinite group $G$
there exists a unique functor
${\bf NR}_G$ from the category of commutative rings with identity
into itself satisfying
the following conditions {\rm:}

\vskip 1mm
$(\text{\bf NR}_G 1)$
As a set
$${\bf NR}_G(R)=\underset{U\leqslant G,\, \text{\rm open}}{{\prod}'}R.$$

\vskip 1mm
$(\text{\bf NR}_G 2)$
For every ring homomorphism $h:A\to B$
and every ${\bf x}=(x_U)'_{U\leqslant G} \in   {\bf NR}_G(A)$ one has
${\bf NR}_G(h)({\bf x})=(f(x_U))'_{U\leqslant G}.$

\vskip 1mm
$(\text{\bf NR}_G 3)$
For every ring $R$
the map
$$\tilde \varphi_R
=\underset{U\leqslant G,\, \text{\rm open}}{{\prod}'} \tilde\varphi_U^{R},$$
where the map $\tilde\varphi_U^{R}:{\bf NR}_G(R)
\to R$ is given by
\begin{equation*}
(x_V)'_V\mapsto
\underset{ U\lesssim V \leqslant G}{{\sum}'}
\varphi_U(G/V)\,x_V\,,
\end{equation*} is a ring homomorphism.
\end{thm}

Eq. \eqref{operation of necklace--equation 1} says that
addition of ${\bf NR}_G(R)$ is defined componentwise
and its multiplication is given by as follows:
for two sequences ${\bf x}=(x_V)'_{V}$ and ${\bf y}=(y_W)'_{W}$
the $U$-th component of the product is given by
\begin{equation*}
({\bf x}\cdot {\bf y})_U:=\underset{V,W}{{\sum}'}
p_V^W(U)\,\, x_V \cdot y_W,
\end{equation*}
where $p_V^W(U)$ represents the number of $VgW$'s in a system of representations
of the double cosets $VgW \subseteq G$ such that
$Z(g,V,W)$ is conjugate to $U$ in $G$.

\vskip 2mm
\noindent{\bf Remark.}

(a) Since $R\supseteq \mathbb Q$ (or a field of characteristic zero)
as a subring, $R$ has a binomial ring structure (see [\cite{O}, Corollary 2.2]).
Indeed, ${\bf NR}_G(R)$ is the necklace ring ${\bf Nr}_G(R)$
of $G$ over $R$ associated with
this $\ld$-ring structure.

(b) If $R$ contains $\mathbb Q$ as a subring then $\tilde
\varphi_R$ is an isomorphism, and if $R$ is a torsion-free
ring which does not contain $\mathbb Q$ as a subring then it is just injective.
However, $\tilde \varphi_R$ is neither injective nor surjective in general.

\vskip 2mm
Let us investigate $\tilde \varphi_R$ (in the above theorem) in more detail.
Since it is $R$-linear, it can be expressed as a matrix form.
Let $P$ be a partially ordered set
consisting of (representatives of conjugacy classes of)
open subgroups of $G$ with the
partial order $\preccurlyeq$ such that
$$V \preccurlyeq W \Leftrightarrow W \lesssim V\,.$$
Consider (and fix) an enumeration
$\{V_i \,:\, 1\le i ,\, V_i \in P\}$
of $P$ satisfying the condition
$$V_i \preccurlyeq V_j \Longrightarrow i\leqq j\,\,.$$
For this enumeration we define
the matrix $\zeta_{P}$ by
\begin{equation*}
\zeta_{P}(V,W):=
\varphi_W(G/V).
\end{equation*}
Using the fact that $\varphi_W(G/V)=0$ unless $V\preccurlyeq W$
we know that $\zeta_{P}$ is a upper-triangular matrix
with the diagonal elements $(N_G(V_i):V_i)$,
the index of $V_i$ in $N_G(V_i)$.
With this notation
\begin{equation}
\tilde \varphi_R({\bf x})={\zeta_{P}}^t \,\,{\bf x},
\text{ where }{\bf x}=
\begin{pmatrix}
x_{V_1}\\
x_{V_2}\\
\vdots
\end{pmatrix}\,\,.
\end{equation}
As a easy consequences of this expression we obtain that
$\tilde \varphi_R$ is invertible if and only if for every $U\in P$
the index $(N_G(V):V)$ is a unit.
When invertible, the inverse of $\tilde \varphi_R$
can be described as follows:
First, let $P(U)$ be a subset of $P$
consisting of open subgroups $V$ of $G$
such that $V\preccurlyeq U$.
Now we set
$$\mu_{P(U)}={\zeta_{P(U)}^{-1}},$$
where ${\zeta_{P(U)}}$ is the matrix obtained from $\zeta_{P}$ by
restricting the index to $P(U)$. Then, for ${\bf a}\in
R^{\underline {\mathcal O(G)}}$
\begin{align*}
\begin{pmatrix}
\tilde\varphi^{-1}({\bf a})_{G}\\
\vdots\\
\tilde\varphi^{-1}({\bf a})_U
\end{pmatrix}
&=\mu_{P(U)}^{t}\,\,
\begin{pmatrix}
a_{G}\\
\vdots \\
a_{U}
\end{pmatrix}
\end{align*}
That is,
\begin{equation}\label{eqn: inverse}
\tilde\varphi^{-1}({\bf a})_U=\sum_{V\preccurlyeq U}\mu_{P(U)}(V,U)\,\,a_V\,.
\end{equation}

\vskip 2mm
If $G$ is abelian and $R$ is torsion-free, the inverse of
$\tilde\varphi_R$ gets much simpler.
For an open subgroup $U$ of $G$ we let $P(U)$ be a partially ordered set
consisting of open subgroups of $G$ containing $U$ with the
partial order $\preccurlyeq$ such that
$$V \preccurlyeq W \Leftrightarrow W\leqslant V\,.$$
Applying M\"obius inversion formula for arbitrary posets
yields that
$$
\tilde\varphi^{-1}({\bf a})_U=
\frac {1}{(G:U)}\sum_{V\preccurlyeq U}\mu_{P(U)}^{\text{ab}}(V,U)\,\,a_V
$$
for ${\bf a}=(a_U)_{U\leqslant G}
\in R^{\underline {\mathcal O(G)}}$
if ${\bf a}$ belongs to the image of $\tilde\varphi.$
Here $\mu_{P(U)}^{\text{ab}}$ is given by
the inverse of the matrix $\zeta_{P(U)}^{\text{ab}}$ given by
\begin{equation*}
\zeta_{P(U)}^{\text{ab}}(V,W):=
\begin{cases}
1, &\text{ if } V\preccurlyeq W,\\
0, &\text{ otherwise.}
\end{cases}
\end{equation*}
Let $R$ contain $\mathbb Q$ as a subring.
Then Theorem \ref{def of the functor:NECKLACE RING}
implies that
for ${\bf a},\,{\bf b}
\in R^{\underline {\mathcal O(G)}}$
\begin{equation}\label{interesting relation of neck:abelian case}
\tilde\varphi^{-1}({\bf ab})_U
={\underset{W,W'\leqslant G \atop W\cap W'=U}{\sum}}
p_W^{W'}(U)\, \tilde\varphi^{-1}({\bf a})_W \,\tilde\varphi^{-1}({\bf b})_{W'}.
\end{equation}
Especially the identity \eqref{interesting relation of neck:abelian case}
comes to us as an interesting and simple formula when applied to
the group $\hat C.$
\begin{prop}\label{formula among coefficients}
Let
$R$ be a commutative ring containing $\mathbb Q$ as a subring.
Then for
${\bf a}=(a_n)_{n\ge 1},{\bf b}=(b_n)_{n\ge 1}
\in R^{\mathbb N}$
the following equation holds {\rm :}
\begin{equation*}
\tilde M({\bf ab},n)
=\sum_{[i,j]=n}(i,i)
\tilde M({\bf a},i)\, \tilde M({\bf b},j),
\end{equation*}
where
$$
\tilde M({\bf a},n):=
\frac 1n\sum_{d|n}\mu(d)a_n.$$
\end{prop}

We now shall describe $\Delta_G(R)$ for a commutative ring $R$.

\vskip 2mm

{\bf Case 1}
Suppose that $R$ contains $\mathbb Q$
(or a field of characteristic zero)
as a subring.
We define $\Delta_G(R)$ by ${\bf NR}_G(R)$.
Proposition \ref{ringhomo of teichmuller} implies that
the map
\begin{align*}
&\tau_R^G:W_G(R)\to \Delta_G(R)\\
&\alpha\mapsto {\underset{U\leqslant G}
{\sum}'}\,\text{Ind}_U^G(R)
\circ \tau^U(R)(\alpha(U)),
\end{align*}
is a ring isomorphism.
Here, for any profinite group $G$ the map
$\tau^G(R): R \to \Delta_G(R)$ is defined by
\begin{equation*}
r\mapsto  \left(M_G(r,U)\right)'_{U\leqslant G},
\end{equation*}
where $M_G(r,U)$ is given by ${\bf M}_G(r,U)$
associated with the binomial ring structure of $R$
(see \cite{OH}).
More explicitly, by Eq. \eqref{eqn: inverse}
and Proposition \ref{ringhomo of teichmuller}
it can be written as
\begin{equation*}
\tilde\varphi^{-1}({\bf a})_U
=\sum_{V\preccurlyeq U}\mu_{P(U)}(V,U)\,\,r^{(G:V)}\,.
\end{equation*}
For example, if $G=\hat C$,
$$M(r,n):=M_{\hat C}(r,\hat C^n)
=\dfrac 1n \sum_{d|n}\mu(d)r^{\frac nd},\,\,n\ge
1.$$

We can also obtain inductions and $\tilde \varphi_R$
as in the previous section.

\vskip 2mm

{\bf Case 2}
Suppose that $R$ does not contain $\mathbb Q$
(or a field of characteristic zero)
as a subring, but it is torsion-free.
We denote by $R\mathbb Q$ the rationalization of $R$, i.e.,
$R\mathbb Q:=R \otimes_{\Z}\mathbb Q.$
Under the natural embedding of $R$ into $R\mathbb Q$
we obtain a map
$$\tau_U({R\mathbb Q})|_R: R \to \Delta_G(R\mathbb Q).$$
By Proposition \ref{multiplicity of teich map} this map is multiplicative.
Combining $\tau_U({R\mathbb Q})|_R$ with the induction map
$$\text{Ind}_U^G(R\mathbb Q): \Delta_U(R\mathbb Q)\to \Delta_G(R\mathbb Q),$$
we obtain a bijective map
$$\tau_{R\mathbb Q}^G|_{W_G(R)}:
W_G(R) \to \text{Im}(\tau_{R\mathbb Q}^G|_{W_G(R)})
$$
defined by
$$\underset {U\leq G}{{\sum}'}\text{Ind}_U^G(R\mathbb Q)
\circ \tau_{R\mathbb Q}^U|_R\,.$$
To be consistent with the notation in \cite {VW}
we shall adopt the following notation
\begin{align*}
&\Delta_G(R):=\left \{{\underset{U}{\sum}'}\,\text{Ind}_U^G(R\mathbb Q)
\circ \tau^U({R\mathbb Q})|_R(\alpha(U))\,:\, \alpha \in W_G(R) \right \}
\subset \Delta_G(R\mathbb Q).
\end{align*}
From that $\tau_{R\mathbb Q}^G$ is a ring isomorphism
and $$\Delta_G(R)=\tau_{R\mathbb Q}^G(W_G(R))$$
it follows the following result.
\begin{prop}
Suppose that $R$ does not contain $\mathbb Q$
{\rm (}or a field of characteristic zero{\rm )}
as a subring, but it is torsion-free.
Then $\Delta_G(R)$ is a subring of $\Delta_G(R\mathbb Q)$,
and moreover it is isomorphic to $W_G(R)$.
\end{prop}

\vskip 2mm
\noindent{\bf Remark.}
Very often it is very important
to know how $\Delta_G(R)$ looks explicitly.
Fortunately this question can be easily answered in some cases.
More precisely, from the argument in [Section 2, \cite{OH}]
it follows that if $R$ is a torsion-free commutative ring
with unity such that $a^p=a$ mod $pR$ if $p$ is a prime
then $R$ has a unique special $\ld$-ring structure
with $\Psi^n=id$ for all $n\ge 1$.
Hence, in this case, Proposition \ref{ringhomo of teichmuller}
implies that
\begin{equation}\label{R:BINOMIAL RING}
\Delta_G(R)=
{\bf Nr}_G(R).
\end{equation}
For example, $\mathbb Z$ or $\mathbb Z_{(r)}$ enjoys this property,
where $\mathbb Z_{(r)}$ means the ring of integers localized at $r$,
that is, $\{m/n \in \mathbb Q\,:\, (n,r)=1\}.$

\vskip 2mm
Next let us discuss exponential maps and inductions.
\begin{lem}\label{welldefinedness of ind}
Suppose that $R$ does not contain $\mathbb Q$
{\rm (}or a field of characteristic zero{\rm )}
as a subring, but it is torsion-free.
Then

{\rm (a)} $\text{\rm Im}(\tau^U(R\mathbb Q)|_R)\subset \Delta_U(R).$

{\rm (b)}
$\text{\rm Im}(\text{\rm Ind}_U^G(R\mathbb Q)|_{\Delta_U(R)})\subset \Delta_G(R).$
\end{lem}

\noindent{\bf Proof.}
(a) is trivial by definition of $\Delta_U(R)$.

(b) By [Section 3.2, \cite{OH}],
for ${\bf a}=(a_W)'_{W\leqslant U}\in W_U(R\mathbb Q),$
\begin{align*}
& (\tau_{R\mathbb Q}^G)^{-1}\circ
\text{Ind}_U^G(R\mathbb Q)\circ \tau_{R\mathbb Q}^U(\bf a)\\
&=(p_V)'_{V\leqslant G}\,,
\end{align*}
where $p_V$ is a polynomial with integral coefficients in those
$a_W $ ($W$ an open subgroup $U$ to which
$W$ is sub-conjugate in $G$).
In view of the definition of $\Delta_G(R),$
this implies that if ${\bf a}=(a_W)'_{W\leqslant U}\in W_U(R)$
then $\tau_{R\mathbb Q}^G((p_V)'_{V\leqslant G})$
should be in $\Delta_G(R)$.
Thus we complete the proof.
\qed

\vskip 2mm
By virtue of Lemma \ref{welldefinedness of ind}
we are able to get the following maps
\begin{equation}\label{torsion-free:exponential map}
\tau^U(R): R\to \Delta_U(R)
\end{equation}
and
$$\text{Ind}_U^G(R): \Delta_U(R)\to \Delta_G(R)$$
from the map $\tau^U(R\mathbb Q)|_R$ and
$\text{Ind}_U^G(R\mathbb Q)|_{\Delta_U(R)}$.
\begin{prop}\label{formula of necklace polynomial}
Suppose that $R$ is torsion-free.
In $R\mathbb Q$,
for an open subgroup $V$ of a profinite group $G$,
and $r, s \in R$,
we have
\begin{equation}\label{rel of neck poly}
M_G(rs,V)
=\underset{W,W'}{{\sum}'}p_W^{W'}(V)\,M_G(r,W)
M_G(s,W')\,.
\end{equation}
If $G$ is abelian, the formula \eqref{rel of neck poly} reduces to
\begin{equation}\label{rel of neck poly: abelian case}
M_G(rs,V)
=\sum_{W,W'\atop
W\cap W'=V}(G:W+W')\,M_G(r,W)M_G(s,W').
\end{equation}
In particular, if $G=\hat C$,
for every positive integer $n$ and $r,s \in R$,
we have
\begin{equation}\label{rel of neck poly: cyclic}
M(rs, n)
=\sum_{[i,j]=n}(i,j)M(r,i) M(s,j)\,,
\end{equation}
where $[i,j]$ is the least common multiple
and $(i,j)$ the greatest common divisor of $i$ and $j.$
\end{prop}
\noindent{\bf Proof.}
From the multiplicativity of $\tau^G(R)$
formula \eqref{rel of neck poly}
follows.
In case $G$ is abelian, the number of all elements in
a system of representations
of the double cosets $WgW'\subseteq G$ where
is given by $(G:W+W')$.
Thus  we have formula \eqref{rel of neck poly: abelian case}.
Finally, if $G=\hat C$,
$$i\mathbb Z \cap j \mathbb Z=[i,j]\mathbb Z,
\quad
i\mathbb Z + j \mathbb Z=(i,j)\mathbb Z.$$
This yields formula \eqref{rel of neck poly: cyclic}.
\qed
\vskip 2mm
Setting $\tau_R^G$ by
\begin{equation}\label{torsion-free:teichmuller}
\underset {U\leq G}{{\sum}'}\text{Ind}_U^G(R)\circ \tau^U(R):
W_G(R)\to \Delta_G(R),
\end{equation}
clearly it is a ring isomorphism.

Finally, from
$\tilde \varphi_{R\mathbb Q}|_{\Delta_G(R)}$
let us obtain the following map
$$\tilde \varphi_{R}
:\Delta_G(R)\to R^{\underline {\mathcal O(G)}}\,.$$
Indeed this map is well-defined since
$\text{Im}(\tilde \varphi_{R\mathbb Q}|_{\Delta_G(R)})
\subset R^{\underline {\mathcal O(G)}}\,\,,$
which follows from the fact
$$\text{Im}\Phi_{R}=\text{Im} \Phi_{R\mathbb Q}|_{W_G(R)}
\text{ and }
\Phi_{R\mathbb Q}
=\tilde \varphi_{R\mathbb Q}\circ\tau_{R\mathbb Q}^G.$$
Consequently we arrive at the following commutative diagram

\begin{equation*}
\begin{picture}(360,70)
\put(100,60){$W_G(R)$}
\put(135,63){\vector(1,0){70}}
\put(210,60){${\Delta}_G(R) \subset {\Delta}_G(R\mathbb Q) $}
\put(160,0){$R^{\underline {\mathcal O(G)}}
\subset {R\mathbb Q}^{\underline {\mathcal O(G)}}$}
\put(215,55){\vector(-1,-1){43}}
\put(120,55){\vector(1,-1) {43}}
\put(165,67){$\tau_R^G$}
\put(165,54){$\cong$}
\put(204,35){$\tilde \varphi_R $}
\put(122,34){$\Phi_R$}
\put(164,38){$\curvearrowright$}
\put(275,55){\vector(-1,-1){43}}
\put(263,33){$\tilde\varphi_{R\mathbb Q} $}
\end{picture}
\end{equation*}

Note that
$\Phi_{R\mathbb Q}, \,\tilde\varphi_{R\mathbb Q}$
are ring-isomorphisms, whereas
$\Phi_{R}, \,\tilde\varphi_{R}$ are injections but not surjections in general.

\vskip 2mm
{\bf Case 3}
Finally we suppose that
$R$ is not torsion-free.
In this case we start by choosing
a surjective ring homomorphism $\rho: R' \to R$
from a torsion free ring $R'$.
For example, we may take $R'=\mathbb Z[R]$ and
$\rho: \mathbb Z[R]\to R$ defined by
\begin{equation}
\sum n_r\cdot e_r \mapsto \sum n_r r,
\end{equation}
where $e_r$ is the basis element of $\mathbb Z[R]$ corresponding to $r\in R.$
By the functoriality of $W_G$ there exists a surjective ring homomorphism
$W_G(\rho) : W_G(R') \to W_G(R)$ defined by
$$(a_U)'_{U\leqslant G}\mapsto (\rho(a_U))'_{U\leqslant G}.$$
Thus it holds
$$\ker W_G(\rho)=W_G(\ker \rho),$$
which implies that the induced map
$$\bar \tau_{R'}^G:W_G(R')/W_G(\ker \rho)
\to \Delta_G(R')/\Delta_G(\ker \rho)$$
from $\tau_{R'}^G$
is an ring isomorphism.
Set
$$\Delta_G(R):= \Delta_G(R')/\Delta_G(\ker \rho)$$
and
$$\tau_R^G:=\bar \tau_{R'}^G.$$
By definition,
$\Delta_G(R)$ and $\tau_R^G$ are well-defined, i.e.,
they are independent
of the choices of $(R',\rho)$ (up to isomorphism).

Let us define exponential maps and inductions.
For an open subgroup $U$ of $G$ it is clear that
the induced map
$$\bar \tau^U(R'): R'
/\ker \rho \to \Delta_U(R')/\Delta_U(\ker \rho)$$
from the map $\tau^U(R')$
is well-defined since $\tau^U(R')(\ker\rho) \subset \Delta_U(\ker\rho).$
Setting $\tau^U(R)$ to be $\bar \tau^U(R')$,
$\tau^U(R)$ can be regarded as a map from
$R$ to $\Delta_U(R).$
In the same manner we can get inductions
$$ \text{Ind}_U^G(R):\Delta_U(R)\to \Delta_G(R).$$
By construction
$\tau^U(R)$ is multiplicative, $\text{Ind}_U^G$ is additive,
and
\begin{equation*}
\tau_R^G=\underset {U\leq G}{{\sum}'}\text{Ind}_U^G(R)\circ \tau^U(R)\,.
\end{equation*}
We can also define a map
$$\tilde \varphi_{R}
:\Delta_G(R')/\Delta_G(\ker \rho) \to
{R'}^{\underline {\mathcal O(G)}}/{(\ker \rho)}^{\underline {\mathcal O(G)}}$$
from $\tilde \varphi_{R'}:\Delta_G(R')\to {R'}^{\underline {\mathcal O(G)}}.$
Regarding $\tilde \varphi_{R}$ as a map from $\Delta_G(R)$
to ${R}^{\underline {\mathcal O(G)}}$, clearly the diagram

\begin{equation}
\begin{picture}(360,70)
\put(100,60){$W_G(R)$}
\put(135,63){\vector(1,0){70}}
\put(210,60){$\Delta_G(R)$}
\put(155,0){$R^{\underline {\mathcal O(G)}}$}
\put(215,55){\vector(-1,-1){43}}
\put(120,55){\vector(1,-1) {43}}
\put(165,67){$\tau_R^G$}
\put(165,54){$\cong$}
\put(204,35){$\tilde\varphi_R$}
\put(122,34){$\Phi_R$}
\put(163,38){$\curvearrowright$}
\end{picture}
\end{equation}
is commutative.
Be careful that in this case
$\Phi_{R}$ and $\tilde\varphi_{R}$ is neither injective
nor surjective in general.

\vskip 5mm
Until now we have described $\Delta_G(R)$ for a given commutative ring $R$.
We shall now describe morphisms.
Given a ring homomorphism $f:A \to B$,
we can define a natural ring homomorphism
\begin{equation*}
\Delta_G(f):  \Delta_G(A) \to  \Delta_G(B)
\end{equation*}
via $\tau_R^G$, i.e., by
\begin{align}\label{def of morphism:delta 0}
\underset {U\leq G}{{\sum}'}\text{Ind}_U^G(A)\circ \tau^U(A)(x_U)
\mapsto \underset {U\leq G}{{\sum}'}\text{Ind}_U^G(B)\circ
\tau^U(B)(f(x_U))
\end{align}
for ${\bf x}=(x_U)'_{U\leqslant G}\in W_G(A).$
Actually if $A$ and $B$ are torsion free, then
\begin{equation}\label{def of morphism:delta}
\Delta_G(f)({\bf x})=(f(x_U))'_{U\leqslant G},
\end{equation}
which can be verified as follows:
First recall in \cite {OH} that
\begin{equation}\label {def of teichmuller*}
\tau_R^G(\alpha)=
\left(\underset {U\leqslant G }{{\sum}'}
\underset {V\leqslant U}{{\sum}'}
M_U(\alpha(U),V)\right)'_{W\leqslant G}\,,
\end{equation}
where $V$ ranges over open subgroups which are conjugate to $W$ in $G$.
Applying the definition \eqref{def of morphism:delta}, we have
\begin{equation*}
f\left(\underset {U\leqslant G }{{\sum}'}
\underset {V\leqslant U}{{\sum}'}
M_U(\alpha(U),V)\right)
=\underset {U\leqslant G }{{\sum}'}
\underset {V\leqslant U}{{\sum}'}
M_U(f(\alpha(U)),V)
\end{equation*}
since $f$ can be viewed as a ring homomorphism from $A\mathbb Q$
to $B\mathbb Q$.
This says that the definition \eqref{def of morphism:delta 0}
coincides with the definition \eqref{def of morphism:delta}.
If $A$ or $B$ is not torsion free,
choose surjective homomorphisms $\rho: A'\to A$ and $\rho':B' \to B$
for torsion free rings $A'$ and $B'$ respectively.
With this situation we have
\begin{lem}\label{def of mor-explicit}
For ${\bf x}\in \Delta_G(A')$
\begin{equation}\label{def of morphism:delta:torsion}
\Delta_G(f)({\bf x}+\Delta_G(\ker \rho))
={\bf y}+ \Delta_G(\ker \rho'),
\end{equation}
where ${\bf y}\in \Delta_G(B')$
is any element satisfying
$$f(x_V+\ker \rho)=y_V+\ker \rho'$$
for every open subgroup $V$ of $G$.
\end{lem}

\noindent{\bf Proof.}
Write
$${\bf x}+\Delta_G(\ker \rho)
=\tau_{A'}^G({\bf z})+\Delta_G(\ker \rho)$$
for some ${\bf z}=(z_U)'_U\in W_G(A').$
\begin{align*}
&\tau_{B'}^G\circ W_G(f)((z_U+\ker \rho)_{U\leqslant G})\\
&=\tau_{B'}^G((f(z_U+\ker \rho))_{U\leqslant G})\\
&=\tau_{B'}^G({\bf z'})+\Delta_G(\ker \rho'),
\end{align*}
where ${\bf z'}=(z'_U)'_U\in W_G(B')$ is any element satisfying
$${\bf z'}+\Delta_G(\ker\rho)=(f(z_V+\ker \rho))_{V\leqslant G}.$$
So we may take
$$y_V=
\underset {U\leqslant G }{{\sum}'}
\underset {V'\leqslant U}{{\sum}'}
M_U(z'_U,V'),
$$
for every open subgroup $V$ of $G$.
Here $V'$ ranges over open subgroups which are conjugate to $V$ in $G$.
Thus we have
\begin{align*}
f(x_V+\ker \rho)
&=f\left(
\underset {U\leqslant G }{{\sum}'}
\underset {V'\leqslant U}{{\sum}'}
M_U(z_U,V')+\ker\rho\right)\\
&=\underset {U\leqslant G }{{\sum}'}
\underset {V'\leqslant U}{{\sum}'}
M_U(f(z_U+\ker \rho),V')\\
&= \underset {U\leqslant G }{{\sum}'}
\underset {V'\leqslant U}{{\sum}'}
M_U(z'_U+\ker \rho,V'))\\
&= \underset {U\leqslant G }{{\sum}'}
\underset {V'\leqslant U}{{\sum}'}
M_U(z'_U,V')+\ker \rho\\
&= y_V+\ker \rho.
\end{align*}
\qed

\vskip 2mm
By definition of $\Delta_G(f)$, the following lemma is
straightforward.
\begin{lem}
$\Delta_G(f)=\tau_B^G\circ W_G(f)\circ {\tau_A^G}^{-1}.$
\end{lem}
As a result, we get the commutativity of the following diagram
\begin{equation}\label{equivalent diagram:Burnside}
\begin{CD}
W_G(A) @>\tau_A^G>\cong>\Delta_G(A)\\
@V{W_G(f)}VV  @V{\Delta_G(f)}VV  \\
W_G(B) @>{\tau_B^G}>\cong> \Delta_G(B)\,\,,
\end{CD}
\end{equation}
which implies that
\begin{thm}\label{summary:Burnside}
The functor $\Delta_G$ is equivalent to the functor $W_G$.
\end{thm}

\noindent{\bf Remark.}
In \cite{Br} Brun gave a new realization of $W_G(R)$
for arbitrary commutative rings with identity and every finite group
using the theory of $G$-Tambara functors.
Briefly speaking it says that the functor $W_G$ coincides with the left adjoint
of the algebraic functor from the category of $G$-Tambara functors to
the category of commutative rings with an action of $G.$
For this purpose he introduced the isomorphism $t_G^R : W_G(R)\to L_GR(G/G)$
given by
\begin{align*}
{(x_U)'}_{U \leq G} \mapsto
\underset{U \leq G}{{\sum}'}L_GR_+(\pi_U^G)\circ
L_GR_{\bullet}(\pi_e^U)(x_U).
\end{align*}
Following the exactly same way as used in [\cite {O}, Section 4],
we can show that essentially
the maps $\tau_R^G,\tau^U(R), \text{Ind}_U^G$ in our notation
coincide with the maps
$$t_G^R, L_GR_{\bullet}(\pi_e^U), L_GR_+(\pi_U^G)$$ respectively
(in Brun's notation).

\vskip 2mm
\section{Aperiodic ring functor ${\bf {\mathcal {AP}}}_G$}

Let $G$ be any profinite group.
In this section we shall construct another
covariant functor ${\bf {\mathcal {AP}}}_G$,
which also will turn out to b equivalent to $W_G$.
Given a commutative ring $R$ with identity, we call
${\bf {\mathcal {AP}}}_G(R)$ {\it the aperiodic ring of $G$ over $R$}.
As in the previous section we shall define
${\bf {\mathcal {AP}}}_G(R)$ differently
according to the three cases.

\vskip 2mm
{\bf Case 1}
Suppose that $R$ is an arbitrary commutative ring
containing $\mathbb Q$ (or a field of characteristic zero)
as a subring.
In this case ${\bf {\mathcal {AP}}}_G(R)$ is defined as follows:

(${\bf {\mathcal {AP}}}_G$1) As a set, it is
$$\underset{U\leqslant G,\, \text{open}}{{\prod}'}R.$$

(${\bf {\mathcal {AP}}}_G$2) Addition is defined componentwise.

(${\bf {\mathcal {AP}}}_G$3) Multiplication is defined so that
for ${\bf x}=(x_V)'_{V}$ and ${\bf y}=(y_W)'_{W}$ in ${\bf {\mathcal {AP}}}_G(R),$
the $U$-th component of ${\bf x} \cdot {\bf y}$ is given by
\begin{equation}\label{multiplication of aperiodic ring}
({\bf x}\cdot {\bf y})_U
:=\underset{V,W \leqslant G}{{\sum}'}\sum_{VgW\subseteq G \atop Z(g,V,W)}
c_U^V(g)x_V y_W,
\end{equation}
where $g$ runs through a system of representations
of the double cosets $UgV \subseteq G$, and
$Z(g,V,W)$ runs over $G$-conjugates to $U,$
and the coefficient $c_U^V(g)$ is given by
$$c_U^V(g)=\frac{(G:Z(g,U,V))}{(G:U)(G:V)}\,\,.$$

\vskip 2mm
\noindent{\bf Remark.}
Let $G$ be an abelian profinite group.
Then the condition that $R$
contains $\mathbb Q$ (or a field of characteristic zero)
as a subring is redundant since
in the multiplication \eqref{multiplication of aperiodic ring}
\begin{equation}\label{abelian-coeff}
\sum_{VgW\subseteq G \atop Z(g,V,W)}
c_U^V(g)
=\frac{(G:U\cap V)(G:U+V)}{(G:U)(G:V)}
=1.
\end{equation}
For example, let $G=\hat C.$
Then the multiplication of ${\bf {\mathcal {AP}}}_G(R)$ is given by
$${\bf (ab)}_n
=\sum_{[i,j]=n}
{a}_i {b}_j
\text{ for } {\bf a,b}\in {\bf {\mathcal {AP}}}_G(R).
$$
Therefore we can form a ring satisfying the above three
conditions for every commutative ring.
Let us denote by ${\bf Ap}_G(R)$ the ring obtained from $R$ in this way.
But in general this is impossible since
coefficients appearing in \eqref{multiplication of aperiodic ring}
do not always have integral values.

\vskip 2mm
Rather than verify that
${\bf {\mathcal {AP}}}_G(R)$ is indeed a commutative ring directly,
we shall prove it indirectly.
To do this, for every conjugacy class of open subgroups $U$ of $G$,
let us define a map
$\varphi^R_U: {\bf {\mathcal {AP}}}_G(R)\to R$ by
$${\bf x} \mapsto
\underset{U\lesssim V \leqslant G}{{\sum}'}
\dfrac{1}{(G:V)}\,\varphi_U(G/V)x_V
$$
for ${\bf x}=(x_U)'_U \, \in  {\bf {\mathcal {AP}}}_G(R)$.
Now we set
$$\varphi_R=\underset{U}{{\prod}'}\varphi^R_U :
\Delta_G(R)\to R^{\underline {\mathcal O(G)}}$$ by $\varphi_R({\bf
x})(U)=\varphi^R_U({\bf x})$ for all ${\bf x} \in \Delta_G(R).$ Note
that if $G$ is abelian, then
\begin{equation}\label{abelian case}
{\varphi_R({\bf x})}(U) = \underset{U\leqslant V \leqslant G}{{\sum}}
x_V
\end{equation}
since $(G:V)=\varphi_U(G/V).$

\begin{prop}{\rm (cf. Varadarajan and Wehrhahn \cite {VW})}
\label{proof of appriodic ring}
For every commutative ring $R$
containing $\mathbb Q$ {\rm (}or a field of characteristic zero{\rm )}
as a subring, we have

{\rm (a)} $\varphi_R$ is an isomorphism of
the additive group ${\bf {\mathcal {AP}}}_G(R)$
onto $R^{\underline {\mathcal O(G)}}$.

{\rm (b)}
For any ${\bf x,\,y} \in {\bf {\mathcal {AP}}}_G(R)$,
$$\varphi_R({\bf xy})=\varphi_R({\bf x})\varphi_R({\bf y})\,.$$
\end{prop}

\noindent{\bf Proof.}
(a) It is clear that $\varphi_R$ is a homomorphism of additive groups.
Let ${\bf x}=(x_U)'_U \, \in  {\bf {\mathcal {AP}}}_G(R)$
satisfy $\varphi_R({\bf x})={\bf 0}$, where ${\bf 0}$ is
the zero element of
$R^{\underline {\mathcal O(G)}}$.
Then for every conjugacy class of open subgroups $U$ of $G$,
$\varphi^R_U({\bf x})=0.$
If $U=G$, then $x_G=0.$
Now assume that $x_V=0$ for all $V$ such that $(G:V)<(G:U).$
Then, from
$$\underset{U\lesssim V \leqslant G}{{\sum}'}
\dfrac{1}{(G:V)}\,\varphi_U(G/V)x_V=0
$$
we can get $x_U=0$.
Thus ${\bf x}={\bf 0}\in {\bf {\mathcal {AP}}}_G(R),$
and which says that $\ker \varphi_R=0.$

Next we will show that $\varphi_R$ is surjective.
For any ${\bf a}=(a_U)'_U \, \in R^{\underline {\mathcal O(G)}}$
we want to find an element ${\bf x}=(x_U)'_U \, \in  {\bf {\mathcal {AP}}}_G(R)$
satisfying
\begin{equation*}
\underset{U\lesssim V \leqslant G}{{\sum}'}
\dfrac{1}{(G:V)}\,\varphi_U(G/V)x_V
=a_U
\end{equation*}
for every conjugacy class of open subgroups $U$ of $G$.
If $U=G$, then $x_G=a_G$.
Let us use mathematical induction on the index.
Assume that we have found $x_V$ for all $V$'s such that $(G:V)<(G:U).$
Then, since
$$
\dfrac{1}{(G:U)}\,\varphi_U(G/U)x_U
=a_U-\underset{U\lesssim V \leqslant G \atop U\ne V}{{\sum}'}
\dfrac{1}{(G:V)}\,\varphi_U(G/V)x_V\,,
$$
$x_U$ is determined by the assumption.
This completes the proof of (a).

(b) For any ${\bf x}=(x_U)'_U \,,\,{\bf y}=(y_V)'_V \in  {\bf {\mathcal {AP}}}_G(R)$,
\begin{align*}
\varphi^R_U({\bf x})\cdot \varphi^R_U({\bf y})
&=\underset{U\lesssim V \leqslant G}{{\sum}'}\frac{1}{(G:V)}\varphi_U(G/V)x_V
\cdot \underset{U\lesssim W \leqslant G}{{\sum}'}\frac{1}{(G:W)}\varphi_U(G/W)y_W\\
&=\underset{U\lesssim V\leqslant G\atop U\lesssim W \leqslant G}{{\sum}'}
\,\frac{x_V y_W}{(G:V)(G:W)}\varphi_U(G/V)\varphi_U(G/W)\,.
\end{align*}
Let $e_U \in R^{\underline {\mathcal O(G)}}$ be the element given by
$e_U=(\delta_{U,V})'_{V\leqslant G}$.
Then
\begin{align*}
\varphi^R_U({\bf xy})
&=\varphi^R_U \left(\underset{V,W}{{\sum}'}\sum_{VgW \subseteq G}
c_U^V(g)x_V y_W \,\,e_{Z(g,V,W)}\right)\\
&=\underset{V,W}{{\sum}'}\sum_{VgW \subseteq G\atop U \lesssim Z(g,V,W)}
c_U^V(g) \, x_V y_W \,\,\varphi^R_U(e_{Z(g,V,W)})\\
&=\underset{V,W}{{\sum}'} \, x_V y_W
\sum_{VgW \subseteq G \atop U \lesssim Z(g,V,W)}\frac {c_U^V(g)}{(G:Z(g,V,W))}
\varphi_U(G/Z(g,V,W))\\
&=\underset{V,W}{{\sum}'} \, \frac {x_V y_W}{(G:V)(G:W)}
\sum_{VgW \subseteq G \atop U \lesssim Z(g,V,W)}
\varphi_U(G/Z(g,V,W)).
\end{align*}
Applying the fact that $\tilde\varphi^R_U$ is a ring homomorphism
it is immediate that
$$\varphi_U^R({\bf xy})=\varphi_U^R({\bf x})\varphi_U^R({\bf y})\,.$$
\qed

\vskip 2mm
As easy but very important results of
Proposition \ref {proof of appriodic ring} we obtain
the following two corollaries.
\begin{cor}\label{aperiodic ring}
Suppose that $R$ is an arbitrary commutative ring
containing $\mathbb Q$ {\rm (}or a field of characteristic zero{\rm )}
as a subring.
Then ${\bf {\mathcal {AP}}}_G(R)$ is a commutative ring.
\end{cor}

\begin{cor}\label{isom of varphi: abelian case}
If $G$ is abelian, then for every commutative ring $R$,
$$\varphi_R: {\bf Ap}_G(R)\to R^{\underline {\mathcal O(G)}}$$
is a ring isomorphism
{\rm (}see the remark containing the identity \eqref{abelian-coeff}{\rm )}.
\end{cor}

\noindent{\bf Proof.}
Recall that if $G$ is abelian, by the equation \eqref{abelian case}
and the remark containing the identity \eqref{abelian-coeff},
we can define ${\bf Ap}_G(R)$ and $\varphi_R$
for any commutative ring $R$.
Moreover, in the proof of
Proposition \ref{proof of appriodic ring},
if we use the identity \eqref{abelian case},
then no non-integer coefficients appear.
This means that
Proposition \ref{proof of appriodic ring} holds for
${\bf Ap}_G(R)$ for every commutative ring $R$.
\qed

\vskip 2mm
As $\tilde\varphi^{-1}$ does,
$\varphi_R^{-1}:R^{\underline {\mathcal O(G)}} \to {\bf Ap}_G(R)$
has a very simple form if $G$ is abelian.
In detail,
$$
\varphi^{-1}({\bf a})_U=
\sum_{V\preccurlyeq U}\mu_{P(U)}^{\text{ab}}(V,U)\,\,a_V
$$
for ${\bf a}=(a_U)_{U\leqslant G}
\in R^{\underline {\mathcal O(G)}}\,\,.$
Actually this observation provides us many interesting relations
among $a_U,\,b_V$'s when applied to the identity
\begin{equation}\label{interesting relation:abelian case}
\varphi^{-1}({\bf ab})_U
={\underset{W,W'\leqslant G \atop W\cap W'=U}{\sum}}
\, \varphi^{-1}({\bf a})_W \,\varphi^{-1}({\bf b})_{W'}
\end{equation}
for ${\bf a},\,{\bf b}
\in R^{\underline {\mathcal O(G)}}\,$.
For example, taking $G=\hat C,$
we can get an analogue of Proposition \ref{formula among coefficients}
\begin{prop}\label{rel of coef-aperioic}
For every commutative ring $R$,
we have the formula
\begin{equation*}
\tilde S({\bf ab},n)
=\sum_{[i,j]=n}
\tilde S({\bf a},i)\, \tilde S({\bf b},j),
\end{equation*}
where
${\bf a}=(a_n)_{n\ge 1},{\bf b}=(b_n)_{n\ge 1}
\in R^{\mathbb N}$
and
$$\tilde S({\bf a},n):=
\sum_{d|n}\mu(d)a_n.$$
\end{prop}

\vskip 2mm
To show that $\Delta_G(R)$ is isomorphic to ${\bf {\mathcal {AP}}}_G(R)$
let us introduce a map
$$\theta_R^G: \Delta_G(R)\to {\bf {\mathcal {AP}}}_G(R)$$
defined by
\begin{equation*}
(x_U)'_{U\leqslant G} \mapsto ((G:U)x_U)'_{U\leqslant G}
\end{equation*}
for all ${\bf x}=(x_U)'_{U\leqslant G}\,.$

\begin{prop}\label{comm of theta and varphi}
Suppose that $R$ is an arbitrary commutative ring
containing $\mathbb Q$ {\rm (}or a field of characteristic zero{\rm)}
as a subring.
Then the map $\theta_R^G$ is a ring isomorphism.
Moreover the following diagram

\begin{equation}
\begin{picture}(360,70)
\put(100,60){$\Delta_G(R)$}
\put(135,63){\vector(1,0){70}}
\put(210,60){${\bf {\mathcal {AP}}}_G(R)$}
\put(155,0){$R^{\underline {\mathcal O(G)}}$}
\put(215,55){\vector(-1,-1){43}}
\put(120,55){\vector(1,-1) {43}}
\put(165,67){$\theta_R^G$}
\put(165,54){$\cong$}
\put(202,35){$\varphi_R$}
\put(123,34){$\tilde \varphi_R$}
\put(163,38){$\curvearrowright$}
\end{picture}
\end{equation}
is commutative.
\end{prop}

\noindent{\bf Proof.}
Since $\tilde\varphi_R$ and $\varphi_R$ are ring-isomorphisms,
for our purpose it suffices to show
$$\tilde\varphi_R= \varphi_R \circ \theta_R^G.$$
For any ${\bf x}=(x_U)'_{U\leqslant G}\, \in \Delta_G(R),$ we get
\begin{align*}
 \varphi^R_U \circ \theta_R^G({\bf x})
&=\underset{U\lesssim V \leqslant G}{{\sum}'}
\dfrac{1}{(G:V)}\,\varphi_U(G/V) (G:V)x_V\\
&=\underset{U\lesssim V \leqslant G}{{\sum}'}
\,\varphi_U(G/V)x_V\\
&=\tilde \varphi^R_U({\bf x}).
\end{align*}
This completes the proof.
\qed

\vskip 2mm
As $\Delta_G(R)$ does, ${\bf {\mathcal {AP}}}_G(R)$
comes equipped with exponential maps and induction maps.
First let us investigate exponential maps.
An {\it exponential} map $\Upsilon^U(R):
R\to {\bf{\mathcal {AP}}}_U(R)$ is defined by
\begin{equation*}
r\mapsto
\left(S_U(r,V)\right)'_{V\leqslant U}\,,
\end{equation*}
where
$S_U(r,V):=(U:V)M_U(r,V)$
for all open subgroups $V$
of $U.$
\begin{prop}
Suppose that $R$ is an arbitrary commutative ring
containing $\mathbb Q$
{\rm (}or a field of characteristic zero{\rm )}
as a subring.
Then $\Upsilon^U(R)$ is multiplicative.
\end{prop}

\noindent{\bf Proof.}
Observe that
\begin{equation}\label{decom of upsilon*}
\Upsilon^U(R)=\theta_R^U\circ \tau^U(R).
\end{equation}
So, the desired result follows from
the fact that $\theta_{R}^U$ and $\tau^U(R)$
are multiplicative.
\qed
\vskip 2mm

For every conjugacy class of open subgroups $U$ of $G$
inductions
$$\text{\bf Ind}_U^G(R): {\bf {\mathcal {AP}}}_U(R) \to {\bf {\mathcal {AP}}}_G(R)$$
are defined by
\begin{equation*}
(x_V)'_{V\leqslant U}\mapsto  (y_W)'_{W\leqslant G},
\end{equation*}
where $y_W$ is the sum of $(G:U)x_V$'s such that
$V$ is conjugate to $W$ in $G$.
\begin{lem}\label{comm of tau and ind*}
Suppose that $R$ is an arbitrary commutative ring
containing $\mathbb Q$ {\rm (}or a field of characteristic zero{\rm )}
as a subring.
Then
\begin{equation}\label{commutativity of ind and theta}
\theta_R^G \circ \text{\rm Ind}_U^G(R)
=\text{\bf Ind}_U^G(R)\circ \theta_R^U\,.
\end{equation}
\end{lem}

\noindent{\bf Proof.}
For any ${\bf x}=(x_V)'_{V\leqslant U} \, \in \Delta_U(R),$
\begin{equation*}
\theta_R^G \circ \text{Ind}_U^G(R)({\bf x})
=(y_W)'_{W\leqslant G},
\end{equation*}
where $y_W=\sum'_V (G:W)x_V.$
On the other hand,
\begin{align*}
\text{\bf Ind}_U^G(R)\circ \theta_R^U({\bf x})
&=\text{\bf Ind}_U^G(R)\left(((U:V)x_V)'_{V\leqslant U}\right)\\
&=\left (\underset{V}{{\sum}'} (G:U)(U:V)x_V \right )'_{W\leqslant G}\\
&=\left ( \underset{V}{{\sum}'} (G:W)x_V\right)'_{W\leqslant G}\\
&=( y_W)'_{W\leqslant G}.
\end{align*}
This completes the proof.
\qed

\vskip 2mm
Composing $\tau_R^G$ with $\theta_R^G$
we get a ring homomorphism from $W_G(R)$ to ${\bf {\mathcal {AP}}}_G(R)$
$$\gamma_R^G:= \theta_R^G\circ\tau_R^G\,.$$
Note that its explicit form is as follows.
\begin{lem}
$$\gamma_R^G={\underset{U\leqslant G}
{\sum}'}\,\text{\bf Ind}_U^G(R)
\circ \Upsilon^U(R).$$
\end{lem}

\noindent{\bf Proof.}
Note that
\begin{align*}
& \theta_{R}^G\circ
{\underset{U\leqslant G}{\sum}'}\,\text{Ind}_U^G(R)
\circ \tau^U(R)\\
&={\underset{U\leqslant G}{\sum}'}\,\text{\bf Ind}_U^G(R)
\circ \theta_R^U
\circ \tau^U(R)\quad\text{ (by Lemma \ref{comm of tau and ind*})}\\
&={\underset{U\leqslant G}{\sum}'}\,\text{\bf Ind}_U^G(R)
\circ \Upsilon^U(R) \quad \text{ (by \eqref{decom of upsilon*})}.
\end{align*}
\qed
\vskip 3mm

{\bf Case 2}
Suppose that $R$ does not contain $\mathbb Q$
(or a field of characteristic zero)
as a subring, but it is torsion-free.
In this case
we define ${\bf {\mathcal {AP}}}_G(R)$ as follows:

Using inductions $\text{\bf Ind}_U^G(R\mathbb Q)$
and exponential maps $\Upsilon^U(R\mathbb Q)$
for all open subgroups $U$ of $G$ simultaneously,
we define
\begin{equation*}
{\bf {\mathcal {AP}}}_G(R):=\left \{{\underset{U\leqslant G}{\sum}'}\,\text{\bf Ind}_U^G
(R\mathbb Q)\circ \Upsilon^U(R\mathbb Q)(\alpha(U))\,:\, \alpha \in W_G(R) \right \}
\subset {\bf {\mathcal {AP}}}_G(R\mathbb Q).
\end{equation*}
From the fact that
$\gamma_{R\mathbb Q}^G$ is a ring isomorphism and
$${\bf {\mathcal {AP}}}_G(R)=\gamma_{R\mathbb Q}^G(W_G(R))$$
it follows that
\begin{prop}
${\bf {\mathcal {AP}}}_G(R)$ is a subring of ${\bf {\mathcal {AP}}}_G(R\mathbb Q).$
\end{prop}

As we did in the case 1, let us define
$$\theta_R^G: \Delta_G(R)\to {\bf {\mathcal {AP}}}_G(R)$$ by
\begin{equation*}
(a_U)'_{U\leqslant G} \mapsto ((G:U)a_U)'_{U\leqslant G}\,.
\end{equation*}
That is, $\theta_R^G
=\theta_{R\mathbb Q}^G|_{\Delta_G(R)}\,.$
In view of the equations \eqref{decom of upsilon*}
and \eqref{commutativity of ind and theta},
$\theta_R^G$ is clearly well-defined.
That is,
$$\theta_{R\mathbb Q}^G
( \Delta_G(R))\subset {\bf {\mathcal {AP}}}_G(R).$$

\begin{prop}\label{isom of theta and: torsion-free}
If $R$ does not contain $\mathbb Q$
{\rm (}or a field of characteristic zero{\rm )}
as a subring, but it is torsion-free,
the map $\theta_R^G$ is a ring isomorphism.
Moreover the following diagram
\begin{equation*}
\begin{picture}(360,70)
\put(100,60){$\Delta_G(R)$}
\put(135,63){\vector(1,0){70}}
\put(210,60){${\bf {\mathcal {AP}}}_G(R) \subset {\bf {\mathcal {AP}}}_G(R\mathbb Q) $}
\put(160,0){$R^{\underline {\mathcal O(G)}}
\subset {R\mathbb Q}^{\underline {\mathcal O(G)}}$}
\put(215,55){\vector(-1,-1){43}}
\put(120,55){\vector(1,-1) {43}}
\put(165,67){$\theta_R^G$}
\put(165,54){$\cong$}
\put(204,35){$\varphi_{R\mathbb Q}|_{{\bf {\mathcal {AP}}}_G(R)} $}
\put(125,34){$\tilde\varphi_{R}$}
\put(164,38){$\curvearrowright$}
\put(285,55){\vector(-1,-1){43}}
\put(270,33){$\varphi_{R\mathbb Q} $}
\end{picture}
\end{equation*}
is commutative.
\end{prop}

\noindent{\bf Proof.}
In view of the definition of ${\bf {\mathcal {AP}}}_G(R)$
it is clear that $\theta_R^G$
is a ring isomorphism.
Moreover we already know that
$\tilde \varphi_R$ comes from
$\tilde \varphi_{R\mathbb Q}|_{\Delta_G(R)}$
and
$$\tilde \varphi_{R\mathbb Q}
=\varphi_{R\mathbb Q}\circ \theta_{R\mathbb Q}^G$$
(by Proposition \ref{comm of theta and varphi}).
Hence
$$\varphi_{R\mathbb Q}({\bf {\mathcal {AP}}}_G(R))
\subset R^{\underline {\mathcal O(G)}}\,,
$$
and the commutativity is immediate.
\qed
\vskip 2mm
Denote $\varphi_{R\mathbb Q}|_{{\bf {\mathcal {AP}}}_G(R)}$
by $\varphi_R:{\bf {\mathcal {AP}}}_G(R)
\to R^{\underline {\mathcal O(G)}}\,.$
Similarly
define exponential maps and inductions by restriction:
\begin{align*}
&\Upsilon^U(R)=
\theta_{R}^U\circ \tau^U(R)\,,\\
&\text{\bf Ind}_U^G(R):
{\bf {\mathcal {AP}}}_U(R)\to {\bf {\mathcal {AP}}}_G(R)\,.
\end{align*}

\begin{lem}\label{welldefinedness of ind:aperiodic case}
Suppose that $R$ does not contain $\mathbb Q$
{\rm (}or a field of characteristic zero{\rm )}
as a subring, but it is torsion-free.
Then $\text{\bf Ind}_U^G(R)$is well-defined, i.e.,
$$\text{\rm Im}
(\text{\bf Ind}_U^G(R\mathbb Q)|_{{\bf {\mathcal {AP}}}_G(R)})
\subset {\bf {\mathcal {AP}}}_G(R).$$
\end{lem}

\noindent{\bf Proof.}
Note that
\begin{align*}
&\text{\bf Ind}_U^G(R\mathbb Q)\circ
\theta_{R\mathbb Q}^U (\Delta_U(R))\\
&=\text{\bf Ind}_U^G(R\mathbb Q) ({\bf {\mathcal {AP}}}_U(R))
\quad \text{ (by Proposition \ref{isom of theta and: torsion-free})}\\
&=\theta_{R\mathbb Q}^G \circ \text{\rm Ind}_U^G(R\mathbb Q)(\Delta_U(R))
\quad \text{ (by Lemma \ref{comm of tau and ind*})}\\
&\subset \theta_{R\mathbb Q}^G (\Delta_G(R))
\quad \text{ (by Lemma \ref{welldefinedness of ind})}\\
&\subset {\bf {\mathcal {AP}}}_G(R).
\end{align*}
\qed

From the multiplicativity of
$\Upsilon^U(R)$
we have an analogue of Proposition \ref{formula of necklace polynomial}.
\begin{prop}
Assume that $R$ is torsion-free.
In $R\mathbb Q$, for $r, s \in R$
and every open subgroup $V$ of $G$,
we have
\begin{equation}\label{multi of s(r,n)}
S_G(rs,V)
={\underset{W,W'}
{\sum}'}\sum_{WgW'\subseteq G\atop
Z(g,W,W')}c_W^{W'}(g)S_G(r,W)\,S_G(s,W')
\end{equation}
where the sum is over $Z(g,W,W')$'s which are conjugate to $V.$
If $G$ is abelian, then the identity \eqref{multi of s(r,n)}
reduces to
\begin{equation}\label{multi of s(r,n);abelian}
S_G(rs,V)
={\underset{W,W'\leqslant G \atop W\cap W'=V}{\sum}}
S_G(r,W)\,S_G(s,W').
\end{equation}
In particular, if $G=\hat C,$
then \eqref{multi of s(r,n);abelian} reduces to the following simple form
\begin{equation}\label{multi of s(r,n);cyclic}
S(rs,n)
={\underset{[i,j]=n}{\sum}}
S(r,i)\,S(s,j)
\end{equation}
for all $n\in \mathbb N.$
\end{prop}

\noindent{\bf Proof.}
Formula \eqref{multi of s(r,n)}
follows from the multiplicativity of
$\Upsilon^U(R)$.
If $G$ is abelian, by applying \eqref{abelian-coeff}
to \eqref{multi of s(r,n)}
we get formula \eqref{multi of s(r,n);abelian}.
Finally formula \eqref{multi of s(r,n);cyclic}
follows from the fact
$i\mathbb Z\cap j\mathbb Z=[i,j]\mathbb Z.$
\qed

\vskip 3mm

{\bf Case 3}
Finally we suppose that
$R$ is not torsion-free.
For a surjective ring homomorphism $\rho: R' \to R$
from a torsion free ring $R'$ let us consider the following diagram

\begin{equation*}
\begin{CD}
\Delta_G(R') @>\theta_{R'}>\cong>{\bf {\mathcal {AP}}}_G(R')\\
@VVV  @VVV  \\
\Delta_G(R)=\Delta_G(R')/\Delta_G(\ker \rho) @>{\bar \theta_{R'}^G}>\cong>
{\bf {\mathcal {AP}}}_G(R')/{\bf {\mathcal {AP}}}_G(\ker \rho)\,\,.
\end{CD}
\end{equation*}

Let
$${\bf {\mathcal {AP}}}_G(R):
={\bf {\mathcal {AP}}}_G(R')/{\bf {\mathcal {AP}}}_G(\ker \rho),$$
and $$\theta_{R}^G:=\bar \theta^G_{R'}.$$
Since $\Delta_G(R)$ is well defined,
${\bf {\mathcal {AP}}}_G(R)$ is well defined, too.
Finally let us discuss how to obtain
exponential maps and inductions, and the map $\varphi_{R}.$
One way to obtain exponential maps $\Upsilon^U(R)$ is to consider the map
\begin{equation*}
\overline {\Upsilon_U^G(R')}:
R'/\ker \rho
\to
{\bf {\mathcal {AP}}}_G(R')/{\bf {\mathcal {AP}}}_G(\ker \rho)\,,
\end{equation*}
and the other is to let
\begin{equation}\label{upsilon commute theta}
\Upsilon^U(R):=\theta^U_R\circ \tau^U(R).
\end{equation}
Of course they coincide with each other.
As for inductions $\text{\bf Ind}_U^G(R)$,
the map $\text{\bf Ind}_U^G(R')$ induces the map
$$\overline {\text{\bf Ind}}_U^G(R'):
{\bf {\mathcal{AP}}}_U(R')/{\bf {\mathcal{AP}}}_U(\ker \rho)
 \to {\bf {\mathcal {AP}}}_G(R')/{\bf {\mathcal {AP}}}_G(\ker \rho).$$
So, let $$\text{\bf Ind}_U^G(R):=\overline {\text{\bf Ind}}_U^G(R').$$
As a final step, from the map
$\varphi_{R'}:{\bf {\mathcal {AP}}}_G(R')\to {R'}^{\underline {\mathcal O(G)}}\,,$
let us derive a map
$$\overline {\varphi_{R'}}
:{\bf {\mathcal {AP}}}_G(R)/{\bf {\mathcal {AP}}}_G(\ker \rho) \to
{R'}^{\underline {\mathcal O(G)}}/{(\ker \rho)}^{\underline {\mathcal O(G)}}\,.$$
Let $\varphi_R :=\overline {\varphi_{R'}}.$
Regarding $\varphi_{R}$ as a map from ${\bf {\mathcal {AP}}}_G(R)$
to ${R}^{\underline {\mathcal O(G)}}$, the diagram

\begin{equation}
\begin{picture}(360,70)
\put(100,60){$\Delta_G(R)$}
\put(135,63){\vector(1,0){70}}
\put(210,60){${\bf {\mathcal {AP}}}_G(R)$}
\put(155,0){$R^{\underline {\mathcal O(G)}}$}
\put(215,55){\vector(-1,-1){43}}
\put(120,55){\vector(1,-1) {43}}
\put(165,67){$\theta_R^G$}
\put(165,54){$\cong$}
\put(204,35){$\varphi_R$}
\put(123,34){$\tilde \varphi_R$}
\put(163,38){$\curvearrowright$}
\end{picture}
\end{equation}
is commutative.
Note that in this case,
$\varphi_{R}, \,\tilde\varphi_{R}$ are neither injective
nor surjective in general.

\begin{lem}
Let $R$ be any commutative ring with identity.
Then we have
\begin{equation}\label{com of induction and upsilon}
\theta_{R}^G\circ \text{\rm Ind}_U^G(R)
=\text{\bf Ind}_U^G
\circ \theta_{R}^U
\end{equation}
and
\begin{equation}\label{com of induction and upsilon*}
\theta_R^G\circ \tau_R^G={\underset{U\leqslant G}{\sum}'}\,\text{\bf Ind}_U^G
\circ \Upsilon^U(R).
\end{equation}
\end{lem}

\noindent{\bf Proof.}
First, let us prove the identity \eqref{com of induction and upsilon}.
If $R$ is torsion-free,
the identity \eqref{commutativity of ind and theta} justifies the desired result.
Now assume that $R$ is not torsion-free.
Then, for any ${\bf x}=(x_V)'_{V\leqslant U} \, \in \Delta_U(R'),$
\begin{align*}
&\theta_R^G \circ \text{Ind}_U^G(R)({\bf x}+\Delta_U(\ker \rho))\\
&=\theta_R^G (\text{Ind}_U^G(R')({\bf x})+\Delta_G(\ker \rho))\\
&=\theta_{R'}^G\circ \text{Ind}_U^G(R')({\bf x})+{\bf {\mathcal {AP}}}_G(\ker \rho)\,.
\end{align*}
On the other hand
\begin{align*}
&\text{\bf Ind}_U^G(R)\circ \theta_R^U({\bf x}+\Delta_U(\ker \rho))\\
&=\text{\bf Ind}_U^G(R)(\theta_{R'}^U ({\bf x})+{\bf {\mathcal {AP}}}_U(\ker \rho))\\
&=\text{\bf Ind}_U^G(R')\circ \theta_{R'}^U ({\bf x})+{\bf {\mathcal {AP}}}_G(\ker \rho))\,.
\end{align*}
By \eqref{commutativity of ind and theta}
$$\theta_{R'}^G\circ \text{Ind}_U^G(R')({\bf x})
=\text{\bf Ind}_U^G(R')\circ \theta_{R'}^U ({\bf x})\,.
$$
Thus we proved \eqref{com of induction and upsilon}.
On the other hand, \eqref{com of induction and upsilon*} follows
from equations
\eqref {upsilon commute theta} and
\eqref{com of induction and upsilon}.
\qed

\vskip 3mm
Given a ring homomorphism $f:A \to B$,
we can define a ring homomorphism
\begin{equation*}
{\bf {\mathcal {AP}}}_G(f): {\bf {\mathcal {AP}}}_G(A) \to
{\bf {\mathcal {AP}}}_G(B)
\end{equation*}
via $\theta_R^G \circ \tau_R^G$, i.e., by
\begin{align*}
\underset {U\leq G}{{\sum}'}\text{\bf Ind}_U^G(A)\circ \Upsilon^U(A)(x_U)
\mapsto \underset {U\leq G}{{\sum}'}\text{\bf Ind}_U^G(B)\circ
\Upsilon^U(B) (f(x_U))
\end{align*}
for ${\bf x}=(x_U)'_{U\leqslant G}\in W_G(R).$
As we show in Section 3,
if $A$ and $B$ are torsion free, it holds that
\begin{equation}\label{def of morphism:aperiodic}
{\bf {\mathcal {AP}}}_G(f)({\bf x})=(f(x_U))'_{U\leqslant G}.
\end{equation}
If $A$ or $B$ is not torsion free,
choose surjective homomorphisms $\rho: A'\to A$ and $\rho':B' \to B$
for torsion free rings $A'$ and $B'$ respectively.
With this situation we can show, as we did in
Lemma \ref{def of mor-explicit},
\begin{lem}
For ${\bf x}\in {\bf {\mathcal {AP}}}_G(A')$
\begin{equation}\label{def of morphism:delta:torsion}
{\bf {\mathcal {AP}}}_G(f)({\bf x}+{\bf {\mathcal {AP}}}_G(\ker \rho))
={\bf y}+ {\bf {\mathcal {AP}}}_G(\ker \rho'),
\end{equation}
where ${\bf y}\in {\bf {\mathcal {AP}}}_G(B')$
is any element satisfying
$$f(x_V+\ker \rho)=y_V+\ker \rho'$$
for every open subgroup $V$ of $G$.
\end{lem}
By definition of $\Delta_G(f)$, the following lemma is
straightforward.
\begin{lem}
${\bf{\mathcal {AP}}}_G(f)
=\theta_B^G\circ \Delta_G(f)\circ {\theta_A^G}^{-1}.$
\end{lem}
In other words we have the commutativity of the following diagram
\begin{equation}\label{equivalent diagram:aperiodic}
\begin{CD}
\Delta_G(A) @>\Delta_G(f)>>\Delta_G(B)\\
@V{\theta_A^G} VV  @V{\theta_B^G} VV \\
{\bf {\mathcal {AP}}}_G(A) @>{\bf {\mathcal {AP}}}_G(f)>>
{\bf {\mathcal {AP}}}_G(B)\,\,.
\end{CD}
\end{equation}

Let $\gamma_R^G:=\theta_R^G \circ \tau_R^G.$
Combining diagrams \eqref{equivalent diagram:Burnside}
and \eqref{equivalent diagram:aperiodic}
with the diagram

\begin{equation}
\begin{picture}(360,70)
\put(100,60){$W_G(R)$}
\put(135,63){\vector(1,0){70}}
\put(210,60){$\Delta_G(R)$}
\put(155,0){${\bf {\mathcal {AP}}}_G(R)$}
\put(215,55){\vector(-1,-1){43}}
\put(120,55){\vector(1,-1) {43}}
\put(165,67){$\tau_R^G$}
\put(165,54){$\cong$}
\put(205,35){$\theta_R^G$}
\put(119,34){$\gamma_R^G$}
\put(163,38){$\curvearrowright$}
\put(180,30){$\cong$}
\put(145,33){$\cong$}
\end{picture}
\end{equation}
we can provide our main result.
\begin{thm}\label{sumary:Aperiodic}
The functor ${\bf {\mathcal {AP}}}_G$ is equivalent to the functor $\Delta_G$.
Hence, it is equivalent to the functor $W_G$.
\end{thm}

\section{Inductions and restrictions}
The functors $W_G,\,\Delta_G$ and ${\bf {\mathcal {AP}}}_G$
come equipped with two important families of additive homomorphisms:
inductions and restrictions.
We already have defined inductions on
$\Delta_G$ and ${\bf {\mathcal {AP}}}_G$ in Section 3 and 4.
On the other hand,
given a commutative ring $R$ with identity and an open subgroup $U$ of $G$
$W_G$ has an additive homomorphism
$v_U:W_U(R) \to W_G(R)$ (see {\cite{DS2,O,OH}).

\begin{prop} \label{preserving induction maps}
For every commutative ring $R$ with identity, the maps
$\tau_R, \,\theta_R$ and $\gamma_R$ preserve induction maps.
That is, for every open subgroup $U$ of $G$ we have\hfill

{\rm (a)} $\text{\rm Ind}_U^G(R)\circ \tau_R^U= \tau_R^G\circ v_U.$

{\rm (b)} $\text{\bf Ind}_U^G(R)\circ \theta_R^U
= \theta_R^G\circ \text{\rm Ind}_U^G(R).$

{\rm (b)} $\text{\bf Ind}_U^G(R)\circ \gamma_R^U
= \gamma_R^G\circ v_U.$
\end{prop}

\noindent{\bf Proof.} We may assume that $R$ contains $\mathbb Q$
as a subring since $\Delta_G(R)$, $\text{\rm Ind}_U^G(R)$, and
$\text{\bf Ind}_U^G(R)$ are constructed from $R\mathbb Q$ if $R$
is torsion free and from a torsion free ring $R'$ if $R$ is not
torsion free. So, assume that $R \supset \mathbb Q.$ In this case
(a) was proven in [Section 2.3, \cite{OH}] and (b) was proven in
Lemma \ref{comm of tau and ind*}. On the other hand, (c) can be
obtained by combining (a) with (b) (by definition of
$\gamma_R^G$).
\qed

\vskip 2mm

For restrictions similar results hold.
Associated with $W_G$
it is known that
there exists a ring homomorphism
$f_U:W_G(R) \to W_U(R)$
for every ring $R$ and every open subgroup $U$ of $G$
(see {\cite{DS2,O,OH}).
Now we define restrictions
$\text{\rm Res}_U^G(R):\Delta_G(R)\to \Delta_U(R)$ as follows:

\vskip 2mm
{\bf Case1}
Suppose that $R$ contains $\mathbb Q$
(or a field of characteristic zero)
as a subring.
Define
$$\text{\rm Res}_U^G(R) :\Delta_G(R)\to \Delta_U(R)$$
by
\begin{equation*}
(b_V)'_{V\leqslant G} \mapsto
\left(\underset{V}{{\sum}'}
\sum_g b_V \right)'_{W \leqslant U},
\end{equation*}
where $g$ ranges over a set
of representatives of $U$-orbits of $G/V$
satisfying $Z(g,U,V)$ is conjugate to $W$ in $U$.

\vskip 2mm
{\bf Case2}
Suppose that $R$ does not contain $\mathbb Q$
(or a field of characteristic zero)
as a subring, but it is torsion free.
From the map
$\text{Res}_U^G(R\mathbb Q)|_{\Delta_G(R)}$
let us obtain the map
$$\text{Res}_U^G(R): \Delta_G(R)\to \Delta_U(R).$$
Actually the following lemma implies that
this map is well-defined.
\begin{lem}\label{welldefinedness of res}
Suppose that $R$ does not contain $\mathbb Q$
{\rm (}or a field of characteristic zero{\rm )}
as a subring, but it is torsion-free.
Then
$$\text{\rm Im}(\text{\rm Res}_U^G
(R\mathbb Q)|_{\Delta_G(R)})\subset \Delta_U(R).$$
\end{lem}

\noindent{\bf Proof.}
In essence the proof is same as that of Lemma \ref{welldefinedness of ind} (b).
By [Section 3.2, \cite{OH}], for ${\bf a}=(a_W)'_{W\leqslant U}\in W_G(R\mathbb Q)$
\begin{align*}
& (\tau_{R\mathbb Q}^U)^{-1}\circ
\text{Res}_U^G(R\mathbb Q)\circ \tau_{R\mathbb Q}^G(\bf a)\\
&=(q_V)'_{V\leqslant U}\,,
\end{align*}
where $q_V$ is a polynomial with integral coefficients in those
$a_{W_i}\,\, 1\le i\le k $
(where, $\{W_i\,:\, 1\le i\le k\}$ is a system of subgroups
of $G$ containing a conjugate of $U$).
In view of the definition of $\Delta_U(R),$
this implies that if ${\bf a}=(a_W)'_{W\leqslant U}\in W_G(R)$
then $\tau_{R\mathbb Q}^U((q_V)'_{V\leqslant U})$
should be in $\Delta_U(R)$.
Thus we complete the proof.
\qed

\vskip 2mm
{\bf Case3}
Suppose that $R$ is not torsion free.
For a surjective ring homomorphism $\rho: R' \to R$
from a torsion free ring $R'$,
the map $\text{\rm Res}_U^G(R') :\Delta_G(R')\to \Delta_U(R')$
induces
$$\overline {\text{\rm Res}_U^G(R')} :
\Delta_G(R')/\Delta_G(\ker \rho) \to \Delta_U(R')/\Delta_U(\ker \rho).$$
Let
$$\text{\rm Res}_U^G(R):= \overline {\text{\rm Res}_U^G(R')}.$$

\vskip 3mm
Associated with the functor ${\bf {\mathcal {AP}}}_G$
we define restrictions
$\text{\bf Res}_U^G(R):{\bf {\mathcal {AP}}}_G(R)\to
{\bf {\mathcal {AP}}}_U(R)$ as follows:

\vskip 2mm
{\bf Case1}
Suppose that $R$ contains $\mathbb Q$
(or a field of characteristic zero)
as a subring.
Define
$$\text{\bf Res}_U^G(R) :
{\bf {\mathcal {AP}}}_G(R)\to {\bf {\mathcal {AP}}}_U(R)$$
by
\begin{equation*}
(b_V)'_{V\leqslant G} \mapsto \left(\underset{V}{{\sum}'}
\sum_g \frac {(U:W)}{(G:V)}\,b_V \right)'_{W \leqslant U},
\end{equation*}
where $g$ ranges over a set
of representatives of $U$-orbits of $G/V$
satisfying $Z(g,U,V)$ is conjugate to $W$ in $U$.

\begin{lem}\label{comm of tau and res*}
Suppose that $R$ is an arbitrary commutative ring
containing $\mathbb Q$ {\rm (}or a field of characteristic zero{\rm )}
as a subring.
Then
\begin{equation}\label{commutativity of res and theta}
\theta_R^U \circ \text{\rm Res}_U^G(R)
=\text{\bf Res}_U^G(R)\circ \theta_R^G\,.
\end{equation}
\end{lem}

\noindent{\bf Proof.}
For any ${\bf x}=(x_V)'_{V\leqslant G} \, \in \Delta_G(R),$
\begin{equation*}
\theta_R^U \circ \text{Res}_U^G(R)({\bf x})
=\left((U:W)\underset{V}{{\sum}'}
\sum_g x_V \right)'_{W \leqslant U},
\end{equation*}
where $g$ ranges over a set
of representatives of $U$-orbits of $G/V$
such that $Z(g,U,V)$ is $U$-conjugate to $W.$
On the other hand,
\begin{align*}
\text{\bf Res}_U^G(R)\circ \theta_R^G({\bf x})
&=\text{\bf Res}_U^G(R)\left(((G:V)x_V)'_{V\leqslant G}\right)\\
&=\left(\underset{V}{{\sum}'}
\sum_g \frac {(U:W)}{(G:V)}\,(G:V)x_V \right)'_{W \leqslant U}\\
&=\left((U:W)\underset{V}{{\sum}'}
\sum_g x_V \right)'_{W \leqslant U}\,.
\end{align*}
This completes the proof.
\qed

\vskip 2mm
{\bf Case2}
Suppose that $R$ does not contain $\mathbb Q$
(or a field of characteristic zero)
as a subring, or it is torsion free.
From the map
$\text{\bf Res}_U^G(R\mathbb Q)|_{{\bf {\mathcal {AP}}}_G(R)}$
let us obtain the map
$$\text{\bf Res}_U^G(R): {\bf {\mathcal {AP}}}_G(R)\to {\bf {\mathcal {AP}}}_U(R).$$
The following lemma guarantees
the well-definedness of this map.
\begin{lem}\label{welldefinedness of res-aperiodic}
Suppose that $R$ does not contain $\mathbb Q$
{\rm (}or a field of characteristic zero{\rm )}
as a subring, but it is torsion-free.
Then
$$\text{\rm Im}(\text{\bf Res}_U^G
(R\mathbb Q)|_{\Delta_G(R)})\subset {\bf {\mathcal {AP}}}_U(R).$$
\end{lem}

\noindent{\bf Proof.}
By virtue of Lemma \ref{comm of tau and res*}
we can apply the method of the proof of
Lemma \ref{welldefinedness of ind:aperiodic case}.
\qed

\vskip 2mm
{\bf Case3}
Finally suppose that $R$ is not torsion free.
For a surjective ring homomorphism $\rho: R' \to R$
from a torsion free ring $R'$,
the map $\text{\bf Res}_U^G(R') :
{\bf {\mathcal {AP}}}_G(R')\to {\bf {\mathcal {AP}}}_U(R')$
induces
$$\overline {\text{\bf Res}_U^G(R')} :
{\bf {\mathcal {AP}}}_G(R')/{\bf {\mathcal {AP}}}_G(\ker \rho)
\to {\bf {\mathcal {AP}}}_U(R')/{\bf{\mathcal {AP}}}_U(\ker \rho).$$
Let
$$\text{\bf Res}_U^G(R):= \overline {\text{\bf Res}_U^G(R')}.$$

\begin{prop}\hfill

{\rm (a)} $\text{\rm Res}_U^G(R)\circ \tau_R^G= \tau_R^U\circ f_U.$

{\rm (b)} $\text{\bf Res}_U^G(R)\circ \theta_R^G
= \theta_R^U\circ \text{\rm Res}_U^G(R).$

{\rm (b)} $\text{\bf Res}_U^G(R)\circ \gamma_R^G
= \gamma_R^U\circ f_U.$
\end{prop}

\noindent{\bf Proof.}
The proof is similar with Proposition \ref{preserving induction maps}.
By definition of ${\bf {\mathcal {AP}}}_G$,
$\text{\rm Res}_U^G$, and $\text{\bf Res}_U^G$
we may assume that $R$ contains $\mathbb Q$ as a subring.
In this case,
(a) was proven in [Section 2.3, \cite{OH}]
and (b) was proven Lemma \ref{comm of tau and res*}.
On the other hand,
(c) follows from (a) and (b) (by definition of $\gamma_R^G$).
\qed

\vskip 2mm
Finally we close this section by remarking
inductions and restrictions on $R^{\mathcal O(G)}$
\begin{align*}
&\nu_U:R^{\mathcal O(U)}\to R^{\mathcal O(G)},\\
&\mathcal F_U:R^{\mathcal O(G)}\to R^{\mathcal O(U)}.
\end{align*}
Note that if $R$ is torsion-free,
then $\tilde \varphi$ is one-to-one.
Therefore in this case we can define $\nu_U$ and $\mathcal F_U$
via $\tilde \varphi$ on the image of $\tilde \varphi$, i.e.,
\begin{align*}
&\nu_U:=\tilde \varphi \circ\text{\rm Ind}_U^G\circ \tilde \varphi^{-1},\\
&\mathcal F_U:=\tilde \varphi \circ\text{\rm Res}_U^G\circ \tilde \varphi^{-1}.
\end{align*}
Note that this method is no more valid unless $R$ is torsion-free.
Now let us consider the case $R$ has torsion.
\begin{lem}
\label{on inverse map}
Let $R$ be a commutative ring
containing $\mathbb Q$
{\rm (}or a field of characteristic zero{\rm )}
as a subring.
For any ${\bf b}=(b_V)'_{V\leqslant G} \in R^{\mathcal O(G)}$
let us write
$$\tilde \varphi^{-1}({\bf b})=(a_V)'_{V\leqslant G}.$$
Then for every open subgroup $V$ of $G$, $a_V$ can be written as
\begin{equation}
\frac {1}{\varphi_V(G/V)}
\underset{V\lesssim W \leqslant G}{{\sum}'}
c_W \,b_W
\end{equation}
for some $c_W \in \mathbb Z.$
\end{lem}

\noindent{\bf Proof.}
In order to prove this we shall
use the mathematical induction on index $(G:V)$.
First note that $a_G=b_G.$
Now we assume that the desired assertion holds for all open subgroups
$W$ of $G$ such that $(G:W)<(G:V).$
From
\begin{equation*}
\underset{V\lesssim W \leqslant G}{{\sum}'}
\varphi_V(G/W)\,a_W=b_V
\end{equation*}
we know that
\begin{align*}
a_V
&=\frac {1}{\varphi_V(G/V)}
\left(b_V-
\underset{V\lesssim W \leqslant G\atop V\nsim W}{{\sum}'}
\varphi_V(G/W)\,a_W\right).
\end{align*}
Note that $\varphi_W(G/W)$ divides $\varphi_V(G/W)$
for every $V\lesssim W$ since the group Aut$(G/W)$ is acting
freely on the set of $G$-morphisms from $G/V$ to $G/W$
and the number of elements of this set equals $\varphi_V(G/W)$
(see \cite{DS2}).
Combining these facts with induction hypothesis yields our assertion.
\qed

\vskip 2mm

\begin{prop}\label{on induction}
Let $R$ be a commutative ring containing $\mathbb Q$
{\rm (}or a field of characteristic zero{\rm )}
as a subring.
For any ${\bf b}=(b_V)'_{V\leqslant U} \in R^{\mathcal O(U)}$
we let
\begin{equation*}
\tilde \varphi \circ\text{\rm Ind}_U^G\circ \tilde \varphi^{-1}
({\bf b})=(t_V)'_{V\leqslant G}.
\end{equation*}
Then for every open subgroup $V$ of $U$, $t_V$ is
a polynomial in
\begin{equation*}
\{b_W\,:\,
V\lesssim W \,{\rm(\text{in }G)} \text{ and }W \leqslant U\}
\end{equation*}
with integer coefficients.
Unless $V$ is not an open subgroup of $U$, then $t_V$ is zero.
\end{prop}

\noindent{\bf Proof.}
Writing
$$\tilde \varphi^{-1}
({\bf b})=(a_{V'})'_{V'\leqslant U},$$
then
$$\text{\rm Ind}_U^G\circ \tilde \varphi^{-1}
({\bf b})=(\alpha_V)'_{V\leqslant G},$$
where $\alpha_V$ is the sum of $a_{V'}$'s such that $V'$ is conjugate to $V$ in $G$.
Therefore
\begin{equation}\label{induction of varphi}
t_V=\underset{V\lesssim W \leqslant G}{{\sum}'}
\varphi_V(G/W)\,\alpha_W.
\end{equation}
Now, we claim that $\varphi_W(U/W)$ divides $\varphi_V(G/W)$.
To show this we observe that $\varphi_W(G/W)$ is $(N_G(W): W)$,
the index of $W$ in the normalizer of $W$ in $G$,
and $N_U(W)=N_G(W)\cap U$.
This implies that
$\varphi_W(U/W)$ divides $\varphi_W(G/W)$.
On the other hand, we already knows that
$\varphi_W(G/W)$ divides $\varphi_V(G/W)$.
Hence we can conclude that
$\varphi_W(U/W)$ divides $\varphi_V(G/W)$.
Now apply Lemma \ref{on inverse map} to Eq. \eqref{induction of varphi}
to get our assertion.
In case where $V$ is not an open subgroup of $U$,
$\alpha_W=0$ for all $V\lesssim W$.
Therefore $t_V$ is zero.
\qed
\vskip 2mm

Proposition \ref{on induction} has an amusing consequence that we can define
inductions $\nu_U$ for arbitrary commutative rings
using polynomials $t_V$'s.
For example, in case $G$ is abelian they are given by
$$(b_V)'_{V\leqslant G} \mapsto (c_W)'_{W \leqslant U},$$
where
\begin{equation*}
c_W=
\begin{cases}(G:U)b_W & \text{ if }W \leqslant U,\\
0 &\text{ otherwise. }
\end{cases}
\end{equation*}
By definition of $\nu_U$ it is straightforward that
\begin{equation*}\label{commutance of nu and ind}
\nu_U \circ \tilde \varphi =\tilde \varphi  \circ \text{\rm
Ind}_U^G.
\end{equation*}
\vskip 2mm
Compared with inductions, restrictions can be defined well.
For every commutative ring,
define $\mathcal F_U:R^{\mathcal
O(G)}\to R^{\mathcal O(U)}$ by
$$(b_V)'_{V\leqslant G}\mapsto (c_W)'_{W\leqslant U}
$$
where
$$c_W:=\begin{cases}
b_V & \text{ if }W \text{ is conjugate to }V \text{ in } G\,,\\
0  & \text{ otherwise.}
\end{cases}$$
Indeed if $R$ contains $\mathbb Q$ as a subring,
then we can verify that
\begin{equation*}
\mathcal F_U=\tilde \varphi \circ\text{\rm Res}_U^G\circ \tilde \varphi^{-1}.
\end{equation*}
Hence we can conclude that
\begin{equation*}\label{commutance of fu and Fu}
\mathcal F_U \circ \tilde \varphi =\tilde \varphi  \circ \text{\rm
Res}_U^G
\end{equation*}
(see \cite {O}).

\section{$q$-deformation of the functor $W$
and its equivalent functors}
Let $R$ be a commutative ring with unity.
Surprisingly it will turn out to exist
well-defined $q$-deformations of
$W(R)$, ${\bf Nr}(R)$, and ${\bf {\mathcal {AP}}}(R)$
(the subscript $G=\hat C$ will be omitted)
for any $q\in \mathbb Z$.
In this section we introduce these $q$-deformations,
and show that their constructions will be functorial
and the resulting functors are equivalent.

For a torsion-free ring $R$,
we let $W^F(R)$ be the group of Witt-vectors over $R$
associated with the formal group law $F$
and $\mathcal C(F,R)$ be the group of curves in $F$.
It is well-known that $W^F(R)$ and $\mathcal C(F,R)$
are values of functors from formal group laws over $R$
to groups, and that
the Artin-Hasse exponential map
\begin{align*}
&H^F:W^F(R)\to \mathcal C(F,R)\\
&\alpha \mapsto {\sum_{n\ge 1}}^F \alpha_n t^n
\end{align*}
gives a natural equivalence of functors
(see \cite{H,L}).
On the other hand,
C. Lenart found in \cite{L} that
it is possible to endow $W^F(R)$ with multiplicative structure
for some formal group laws, more precisely for
$$F_q(X,Y):=X+Y-qXY, \quad q\in \mathbb Z.$$
To avoid confusion
we shall adapt and use all notations and definitions in \cite {L}
without changes (and without any explanations).
Fix a formal group law $F_q$ for some $q\in \mathbb Z$,
and define a map $\lambda:{\bf Gh}(R)\to {\bf Gh}(R)$
by $\lambda({\bf x})=(nx_n)_{n\ge 1}$
for all ${\bf x}=(x_n)_{n\ge 1}.$
With this notation, C. Lenart proved the following facts on
the multiplicative structure of $W^q(R)$
(the superscript $q$ will be used instead of $F_q$).

\begin{lem}{\rm (Lenart \cite{L})}\label{lenart}

{\rm (a)}
There is a ring structure on $W^q(\mathbb Z)$ such that
the restriction of $\lambda \circ w^q$ is a ring homomorphism.
The map $H^q$ provides an isomorphic ring structure on
$\mathcal C(F_q,\mathbb Z).$

{\rm (b)}
Let $R$ be a torsion free commutative ring with identity.
There are ring structures on
${\bf Nr}^q(R),\,\,W^q(R)$, and $\mathcal C(F_q,R)$ such that
the restrictions of $\lambda\circ g^q$ and $\lambda\circ w^q$
are ring homomorphism, and the restriction of
$H^q$ is a ring isomorphism.
The corresponding Frobenius and the Verschiebung operators
are endomorphisms of these rings.

{\rm (c)}
If $R$ is a $\mathbb Q$-algebra or $R=\mathbb Z$
or $R=\mathbb Z_{(r)}$, there are restrictions of the maps
$T^q$ and $c^q$ to the rings ${\bf Nr}^q(R),\,\,W^q(R)$,
and $\mathcal C(F_q,R)$, and they are isomorphisms.
Furthermore, the Verschiebung and Frobenius operators commutes with
the maps $T^q$ and $c^q.$

{\rm (d)} The statements in {\rm (b)} also holds
if $q$ is viewed as an indeterminate, and $R$ is an algebra
over the ring of numerical polynomials in $q$
{\rm (}in particular, if it coincides with this ring{\rm )}.
\end{lem}

Actually Lemma \ref{lenart} implies a little more.
Let $A=\mathbb Z[a_1,a_2,\cdots;b_1,b_2,\cdots],$
and let ${\bf a}=(a_n)_{n\ge 1}$
and ${\bf b}=(b_n)_{n\ge 1}.$
Define
\begin{align*}
&\Phi_{A}^q:W^q(A)\to {\bf Gh}(A),\\
&{\bf x}\mapsto
\left(\sum_{d|n}dq^{\frac nd -1}x_d^{\frac nd}\right)_{n\ge 1}\,\,.
\end{align*}
Let us solve the following equations
\begin{align*}
&\Phi_A^q({\bf s}^q)
=\Phi_A^q({\bf a})+
\Phi_A^q({\bf b}),\\
&\Phi_A^q({\bf p}^q)
=\Phi_A^q({\bf a})\cdot
\Phi_A^q({\bf b}),\\
&\Phi_A^q({\bf \iota}^q)
=-\Phi_A^q({\bf a}),
\end{align*}
where ${\bf s}^q=(s_n^q)_{n\ge 1}$,
${\bf p}^q=(p_n^q)_{n\ge 1},$ and ${\bf \iota}^q=(\iota_n^q)_{n\ge 1}.$
Note that on these equations Lemma \ref{lenart} (a) implies the following fact.
\begin{lem}{\rm (cf. Lenart \cite{L})}\label{generalize of witt ring}
Fix $q\in \mathbb Z$.
Then $s_n^q,\,p_n^q$ are polynomials in $a_d,b_d$'s for $d\,|\,n$
with integral coefficients.
Similarly,
$\iota_n^q$ is a polynomial in $a_d$'s for $1\le d\le n$
with integral coefficients for every $n\ge 1$.
\end{lem}
\noindent{\bf Proof.}
The desired result for $p_n^q$ was due to M. Hopkins (see \cite{L}).
So we shall consider $s_n^q$ and $\iota_n^q$ only.
To begin with,
we note that $$\Phi_A^q({\bf x})=\frac 1q \Phi_A(q{\bf x}).$$
Thus we have
\begin{align*}
\Phi_A^q({\bf a})+\Phi_A^q({\bf b})
&=\frac 1q \Phi_A( q{\bf a})+\frac 1q \Phi_A( q{\bf b})& \text{ in }{\bf Gh}(A)\\
&=\frac 1q \Phi_A\left(q{\bf a}+q{\bf b}\right)& \text{ in }W(A)\\
&=\Phi_A^q\left( \frac 1q(q{\bf a}+q{\bf b})\right)
& \text{ in }W(A).
\end{align*}
This computation implies that
$$s_n^q=\frac 1q s_n(qa_d,qb_d \,:\, d|n).$$
Since $s_n$ is a polynomial without constant term
$s_n^q$ is a polynomial with integer constants.
Similarly
\begin{align*}
-\Phi_A^q({\bf a})
&=-\frac 1q \Phi_A( q{\bf a})\\
&=\frac 1q \Phi_A(-q{\bf a}) \quad (`-' \text{means the inverse of $+$ in }W(A))\\
&=\Phi_A^q\left( \frac 1q (-q{\bf a})\right).
\end{align*}
From this it follows that
$$\iota_n^q=\frac 1q \iota_n(-qa_n \,:\, 1\le d \le n).$$
Clearly it has integer coefficients
since $\iota_n$ has integer coefficients and no constant term.
\qed

Lemma \ref{generalize of witt ring} has an amusing consequence that
we can define $W^q$ as a functor from the category
of commutative rings with identity to the category of
commutative rings.

\noindent{\bf Remark.}
In general $W^q(R)$ does not have the identity unless $R$
contains $\mathbb Q$ as a subring.
Indeed, if exists,
the identity can be determined inductively by letting
$$\sum_{d|n}dq^{\frac nd -1}a_d^{\frac nd}=1$$
for all $n.$
\vskip 2mm
For $q\in\mathbb Z$, let us define $W^q(R)$
for any commutative ring $R$ with identity
as follows:

($W^q$1) As a set, it is $R^{\mathbb N}.$

({$W^q$2) For any ring homomorphism  $f:A\to B$, the map
$W^q(f):{\bf a}\mapsto (f(a_n))_{n\ge 1}$ is a ring homomorphism
for ${\bf a}=(a_n)_{n\ge 1}$.

($W^q$3) The map $\Phi_R^q: W^q(R)\to {\bf Gh}(R)$ defined by
$${\bf a} \mapsto
\left(\sum_{d|n}dq^{\frac nd -1}a_d^{\frac nd}\right)_{n\ge 1}
\text{  for }{\bf a}=(a_n)_{n\ge 1}$$
is a ring homomorphism.

On the other hand, since
Artin-Hasse map $H^q$ is determined by universal polynomials
with integer coefficients
\begin{align*}
H^q({\bf x})=x_1t+x_2t+(x_3-qx_1x_2)t^3+(x_4-qx_1x_3)t^4+\cdots\,,
\end{align*}
we can endow $\mathcal C(F_q,R)$ with
the ring structure via $H^q$.
Consequently for every commutative ring $R$ with identity
we get a ring isomorphism $H^q: W^q(R) \to \mathcal C(F_q,R).$
Frobenius and Verschiebung operators
can be also defined $W^q(R)$ and $\mathcal C(F_q,R)$,
which are preserved by $H^q$
since they are also given by universal polynomials
with integer coefficients in torsion free cases.

As we did in Section 3 and 4, we can construct functors $\Delta^q$
and ${\bf {\mathcal {AP}}}^q$
which are equivalent to $W^q$ (so to $\mathcal C(F_q,\, \cdot)$).
Indeed the process of their construction is exactly same as
that of $\Delta_G$ and ${\bf {\mathcal {AP}}}_G$.
First we define $\Delta^q(R)$ as follows:

\vskip 2mm
{\bf Case1}
If $R$ contains $\mathbb Q$ as a subring, then we let
$$\Delta^q(R)={\bf Nr}^q(R).$$
The ring
${\bf Nr}^q(R)$ is defined by the following conditions
(for the completeness, refer to \cite{L}):

(${\bf Nr}^q$1) As a set, it is $R^{\mathbb N}.$

(${\bf Nr}^q$2) Addition is defined componentwise.

(${\bf Nr}^q$3) Multiplication is defined so that
for ${\bf x}=(x_n)_{n\ge 1}$ and ${\bf y}=(y_n)_{n\ge 1}$ in
${\bf Nr}^q(R),$
the $n$-th component of ${\bf x} \cdot {\bf y}$ is given by
$$\sum_{[i,j]|n}(i,j)P_{n,i,j}(q)x_i y_j,$$
where $P_{n,i,j}(q),\,[i,j]\,|\,n$ are numerical polynomials in
$\mathbb Q[q]$ given by
$$\frac{j}{(i,j)q}
\sum_{d\,|\, n/[i,j]}
\tau^q \left(\frac {n}{[i,j]d}\,,\frac {n}{i}\right)
S(q^{[i,j]/j},d),$$
and the notation $\tau^q(i,n)$ denotes the quantity
$\sum_{d|i}\mu^q(1,d)\zeta^q(d,n)$
(for the definition of $\mu^q$ and $\zeta^q$, see below).
Introduce the $q$-exponential map $M^q:R\to {\bf Nr}^q(R)$ defined by
$${\bf x}\mapsto (M^q(x,n))_{n\ge 1},$$
where
$$M^q({\bf x},n)=\sum_{d|n}\mu^q(d,n)\frac {q^{d-1}}{d}x^d.$$
Here $\mu^q(d,n)$ is the $(d,n)$-th entry of the inverse of
the matrix $\zeta^q$ defined on the lattice $D(n)$ of divisors
given by
\begin{equation*}
\zeta^q(d_1,d_2):=
\begin{cases}
\frac {d_1}{d_2}\, q^{\frac {d_2}{d_1}-1} & \text{ if } d_1|d_2,\\
0 & \text{ otherwise.}
\end{cases}
\end{equation*}
Note that unless $q=1$, $M^q$ is not multiplicative.
Denote Frobenius and Verschiebung operators by $V_r^q$
and $f_r^q$ for $r\ge 1$ respectively.
Note that Verschiebung operators are defined regardless of $q$,
that is, $V_r^q=V_r^1=V_r$, whereas $f_r^q$,
which is a ring homomorphism,
is defined by
\begin{equation}\label{frobeni of q-version}
{\bf x}\mapsto \left(r\sum_{d|rn}\tau^q
\left(\frac {rn}{[r,d]},\frac{rn}{d}\right)x_d\right)_{n\ge 1}.
\end{equation}
From \cite {L} it follows that the $q$-Teichm\"uller map
\begin{align*}
&T^q:W^q(R) \to {\bf Nr}^q(R)\\
&{\bf x}\mapsto \sum_{n=1}^{\infty}V_nM^q(x_n),
\quad {\bf x}=(x_n)_{n\ge 1}
\end{align*}
is a ring isomorphism.

Finally we let $\tilde\varphi_R^q=\lambda\circ g^q$.
More precisely
$\tilde\varphi_R^q: {\Delta}^q(R) \to {\bf Gh}(R)$ is defined by
$$
{\bf x} \mapsto
\left(\sum_{d|n}dq^{\frac nd -1}x_d \right)_{n\ge 1}
$$
for ${\bf x}=(x_n)_{n\ge 1} \, \in {\Delta}^q(R)$.
It is well known (see \cite{H,L}) that
$$\Phi_R^q=\tilde\varphi_R^q \circ T^q,$$
and moreover Lemma \ref{lenart} says that
$\tilde\varphi_R^q$ is a ring homomorphism.
\vskip 2mm
{\bf Case2}
If $R$ does not contain $\mathbb Q$ as a subring, but
it is torsion free, then we let
$$\Delta^q(R):=T^q(W^q(R)).$$
By definition $\Delta^q(R\mathbb Q)$
is isomorphic to $W^q(R)$ and
$\mathcal C(F_q,R).$
It is easy to verify that
there exist restrictions of the maps
$V_r:\Delta^q(R\mathbb Q)\to \Delta^q(R\mathbb Q)$,
$f_r^q:\Delta^q(R\mathbb Q)\to \Delta^q(R\mathbb Q)$,
$M^q:R\mathbb Q \to {\bf Nr}^q(R\mathbb Q)$,
and $\tilde\varphi_{R\mathbb Q}^q: {\Delta}^q(R\mathbb Q) \to {\bf Gh}(R\mathbb Q)$
to the maps
$V_r:\Delta^q(R)\to \Delta^q(R)$,
$f_r^q:\Delta^q(R)\to \Delta^q(R)$,
$M^q:R\to {\bf Nr}^q(R)$,
and $\tilde\varphi_R^q: {\Delta}^q(R) \to {\bf Gh}(R)$.
For example, for Frobenius operators observe that
the equation \eqref{frobeni of q-version} has
only integer coefficients, i.e.,
$$r\tau^q
\left( \frac {rn}{[r,d]},\frac{rn}{d}\right)\in \mathbb Z$$
for all $q\in \mathbb Z.$
This implies that
$f_r^q(\Delta^q(R)) \subset \Delta^q(R).$

\begin{prop}
If $R$ is a torsion-free commutative ring
with unity such that $a^p=a$ mod $pR$ if $p$ is a prime,
then
$$
\Delta^q(R)={\bf Nr}^q(R).
$$
\end{prop}

\noindent{\bf Proof.}
By \cite{L} we know that
$$M^q(x,n)
=\sum_{d|n}\left(d\tau^q\left(\frac nd,n\right)\right)M(x,d),$$
and $d\tau^q\left(\frac nd,n\right)\in \mathbb Q[q]$
are numerical polynomials for all positive integers $d,\,n$ with $d|n.$
On the other hand, by hypothesis on $R$,
it has a unique special $\ld$-ring structure
with $\Psi^n=id$ for all $n\ge 1$
(see [Section 2, \cite{OH}]).
It implies that $M(x,d)\in R$ for $x\in R.$
Therefore we get $M^q(x,n)\in R,$ and which means that
$\Delta^q(R)\subset{\bf Nr}^q(R).$
Now, to show $\Delta^q(R)\supset{\bf Nr}^q(R)$
choose an arbitrary element ${\bf x}=(x_n)_{n\ge 1}\in {\bf Nr}^q(R).$
As in the classical case (i.e., $q$=1), we can find
$a_n$'s in $R$ inductively such that $T^q((a_n)_{n\ge 1})={\bf x}.$
This completes the proof.
\qed
\vskip 2mm
In case $R$ is torsion-free
$\tilde\varphi_R^q$ is an injective ring homomorphism and its inverse
is given by
\begin{equation}
(\tilde\varphi_R^q)^{-1}({\bf x})=\sum_{d|n}\mu^q(d,n)\,\frac {x_d}{d}
\end{equation}
if ${\bf x}$ belongs to the image of $\tilde\varphi_R^q$.
Thus we can state a $q$-analogue of
Proposition \ref{formula among coefficients}.
\begin{prop}\label{formula among coefficients-q}
Let $R$ be a torsion-free ring. Then for
${\bf a}=(a_n)_{n\ge 1},{\bf b}=(b_n)_{n\ge 1}
\in \text{\rm Im}\,\tilde \varphi_R^q$
the following equation holds {\rm(}in $R\mathbb Q${\rm)}{\rm :}
\begin{equation*}
\tilde M^q({\bf ab},n)
=\sum_{[i,j]=n}(i,i)
P_{n,i,j}(q)\tilde M^q({\bf a},i)\, \tilde M^q({\bf b},j),
\end{equation*}
where
$$
\tilde M({\bf a},n):=\sum_{d|n}\mu^q(d,n)\frac {a_d}{d}.$$
\end{prop}
In particular, by considering the case
${\bf a}=(q^{n-1}x^n)_{n\ge 1},\,{\bf b}=(q^{n-1}y^n)_{n\ge 1}$,
we can recover Proposition 5.15 in \cite {L}
$$M^q(qxy)=qM^q(x)\cdot M^q(y).$$
\vskip 2mm
{\bf Case3}
Finally if $R$ is not torsion-free,
for a surjective ring homomorphism $\rho: R' \to R$
from a torsion free ring $R'$, we define
\begin{align*}
\Delta^q(R):=\Delta^q(R')/\Delta^q (\ker \rho).
\end{align*}
Let us obtain $T^q:W^q(R)\to \Delta^q(R)$,
$V_r$, $f_r^q$ ($r\ge 1$), $M^q$, and $\tilde \varphi_R^q$ in the same way
as we did in Section 3.

Consequently we arrive at the following commutative diagram

\begin{equation*}\label{commuting diagram :formal group**}
\begin{picture}(360,70)
\put(100,60){$W^q(R)$}
\put(133,63){\vector(1,0){30}}
\put(170,60){${\Delta}^q(R)$}
\put(170,1){${\bf Gh}(R)$}
\put(205,63){\vector(1,0){30}}
\put(200,67){$H^q\circ {T^q}^{-1}$}
\put(247,60){$\mathcal C(F_q,R)$}
\put(185,56){\vector(0,-1){40}}
\put(120,55){\vector(1,-1){45}}
\put(207,12){\vector(1,1){45}}
\put(240,35){$\lambda^{-1}\circ E^q$}
\put(145,67){$T^q$}
\put(145,55){$\cong$}
\put(215,55){$\cong$}
\put(170,35){$\tilde \varphi_R^q$}
\put(120,35){$\Phi_R^q$}
\end{picture}
\end{equation*}

Note that all maps appearing in this diagram preserve each
Frobenius and Verschiebung operator since
they do in case $R=\mathbb Z$.
As for morphisms of
$\Delta^q$, they can be obtained in the
same way as done in Section 3.

Finally we are going to introduce $q$-aperiodic ring functor
${\bf {\mathcal {AP}}}^q$.
To begin with, we define ${\bf Ap}^q (R)$ as follows.

(${\bf Ap}^q$1) As a set, it is $R^{\mathbb N}.$

(${\bf Ap}^q$2) Addition is defined componentwise.

(${\bf Ap}^q$3) Multiplication is defined so that
for ${\bf x}=(x_n)_n$ and ${\bf y}=(y_n)_n$ in
${\bf Ap}^q(R),$
the $n$-th component of ${\bf x} \cdot {\bf y}$ is given by
$$\sum_{[i,j]|n}\frac {n}{[i,j]}P_{n,i,j}(q)x_i y_j.$$

Let us define a map
$\varphi_R^q: {\bf Ap}^q(R) \to {\bf Gh}(R)$ by
$$
{\bf x} \mapsto
\left(\sum_{d|n}q^{\frac nd -1}x_d \right)_{n\ge 1}
$$
for ${\bf x}=(x_n)_{n\ge 1} \, \in {\bf Ap}^q(R)$.

\begin{prop}\label{proof of appriodic ring*}
For every commutative ring $R$ wit identity, we have

{\rm (a)} $\varphi_R^q$ is an isomorphism of
the additive group ${\bf Ap}^q(R)$
onto ${\bf Gh}(R)$.

{\rm (b)}
For any ${\bf x,\,y} \in {\bf Ap}^q(R)$,
$$\varphi_R^q({\bf xy})=\varphi_R^q({\bf x})\varphi_R^q({\bf y})\,.$$
\end{prop}

\noindent{\bf Proof.}
(a) This proof is identical to that of
Proposition \ref{proof of appriodic ring} (a).

(b) Let ${\bf x}=(x_n)_{n\ge 1}$ and ${\bf x}=(x_n)_{n\ge 1}.$
To prove our assertion we have to show that for every $n\ge 1$
it holds
\begin{align*}
\sum_{d|n}q^{\frac nd -1}\left(\sum_{[i,j]|d}P_{d,i,j}(q)x_iy_j \right)
=\left(\sum_{e|n}q^{\frac ne -1}x_e \right)
\left(\sum_{f|n}q^{\frac nf -1}y_f \right).
\end{align*}
Equivalently, we have only to show that for every $n\ge 1$ and $i,j\,|\,n$
\begin{align}\label{ring homo of q-aperiodic ri}
\sum_{[i,j]|d \atop d|n}\frac {d}{[i,j]}q^{\frac nd -1}
P_{d,i,j}(q)=q^{\frac ni +\frac nj -2}.
\end{align}
This follows from the fact that $\tilde\varphi_R^q$ is a ring homomorphism.
In fact, by computing the coefficient of $x_iy_j$ in
$\tilde\varphi_R^q({\bf x\cdot y})$
and $\tilde\varphi_R^q({\bf x})\tilde\varphi_R^q({\bf y}),$
we obtain
\begin{align*}
\sum_{d|n}dq^{\frac nd -1}\left(\sum_{[i,j]|d}(i,j)P_{d,i,j}(q)\right)
=ijq^{\frac ni +\frac nj -2},
\end{align*}
and which is clearly equivalent to \eqref{ring homo of q-aperiodic ri}.
\qed

\vskip 2mm
Note that
the inverse of $\varphi_R^q$
is given by
\begin{equation}
(\varphi_R^q)^{-1}({\bf x})=\sum_{d|n}\mu^q(d,n) \,\frac nd x_d\,.
\end{equation}
Thus we have a $q$-analogue of
Proposition \ref{rel of coef-aperioic}.
\begin{prop}
Let $R$ be a commutative ring.
Then for
${\bf a}=(a_n)_{n\ge 1},{\bf b}=(b_n)_{n\ge 1}
\in R^{\mathbb N}$
the following equation holds {\rm :}
\begin{equation*}
\tilde S^q({\bf ab},n)
=\sum_{[i,j]=n}
P_{n,i,j}(q)\,\tilde S^q({\bf a},i)\, \tilde S^q({\bf b},j),
\end{equation*}
where
$$
\tilde S({\bf a},n):=\sum_{d|n}\mu^q(d,n)\,\frac nd \,a_d.$$
\end{prop}
\vskip 2mm
According to Proposition \ref{proof of appriodic ring*}
${\bf Ap}^q(R)$ is a commutative ring.
As in Section 4, we define
$${\bf V}_r: {\bf Ap}^q(R)\to {\bf Ap}^q(R)$$
by
$$(x_n)_{n\ge 1}\mapsto (rx_{\frac nr})_{n\ge 1}
\qquad
\text{ with } x_{\frac nr}:=0 \text{ if }\frac nr \notin \mathbb N,$$
and define
$${\bf f}_r^q: {\bf Ap}^q(R)\to {\bf Ap}^q(R)$$
by
\begin{equation}\label{frobeni of q-version:aperiodic}
{\bf x}\mapsto \left(r\sum_{d|rn}\tau^q
\left(\frac {rn}{[r,d]},\frac{rn}{d}\right)\frac nd x_d\right)_{n\ge 1}.
\end{equation}
Actually these operators are defined via the isomorphism $\theta_R^q$
(see {\bf case 1} in the below).
This means that $\theta_R^q$ is compatible with these operators.
Now we are ready to construct the functor ${\bf {\mathcal{AP}}}^q$.
Let us define ${\bf {\mathcal{AP}}}^q(R)$ in the following steps.

{\bf Case 1}
Suppose that $R$ is an arbitrary commutative ring
containing $\mathbb Q$ (or a field of characteristic zero)
as a subring.
In this case we let
$${\bf {\mathcal{AP}}}^q (R)={\bf Ap}^q (R).$$
From Proposition \ref{proof of appriodic ring*}
it follows that the map
$$\theta_R^q: \Delta^q(R)\to {\bf {\mathcal{AP}}}^q (R)$$ given by
\begin{equation*}
{\bf x} \mapsto (nx_n)_{n\ge 1}
\end{equation*}
for all ${\bf x}=(x_n)_{n\ge 1}$
is a ring isomorphism, and
$\tilde\varphi_R^q=\varphi_R^q \circ\theta_R^q.$
We already mentioned that $\theta_R^q$ preserves
the Frobenius and Verschiebung operators.
Now, composing $\theta_R^q$ with $T^q$
we get the isomorphism
$$\sum_{n=1}^{\infty}{\bf V}_n\circ S^q: W^q(R)\to {\bf {\mathcal{AP}}}^q (R).$$
Here the $q$-exponential map $S^q:R\to {\bf {\mathcal{AP}}}^q(R)$
is given by
$${\bf x}\mapsto (S^q({\bf x},n))_{n\ge 1},$$
where
$$S^q({\bf x},n):=nM^q({\bf x},n).$$

{\bf Case2}
If $R$ does not contain $\mathbb Q$ as a subring, but
it is torsion free, then we let
$${\bf {\mathcal{AP}}}^q (R):=\theta_{R\mathbb Q}^q(\Delta^q(R)).$$
Let us obtain $\theta_R^q:\Delta^q(R)\to {\bf {\mathcal{AP}}}^q(R)$,
${\bf V}_r$, ${\bf f}_r^q$ ($r\ge 1$), $S^q$,
and $\varphi_R^q$ in the same way via restrictions
as we did in Section 4.
\vskip 2mm
{\bf Case3}
Finally if $R$ is not torsion-free,
for a surjective ring homomorphism $\rho: R' \to R$
from a torsion free ring $R'$, we define
\begin{align*}
{\bf {\mathcal{AP}}}^q (R):
={\bf {\mathcal{AP}}}^q (R')/{\bf {\mathcal{AP}}}^q (\ker\rho).
\end{align*}

\vskip 2mm
By construction $\Delta_q(R)$ is isomorphic to ${\bf {\mathcal{AP}}}^q (R)$
for every commutative ring $R$ with identity,
and the process to obtain the Frobenius
and Verschiebung operators and $\varphi_R^q$
seems to be routine. So we shall skip it.
Similarly it can be shown that the following diagram

\begin{equation}
\begin{picture}(360,70)
\put(100,60){$\Delta_q(R)$}
\put(135,63){\vector(1,0){70}}
\put(210,60){${\bf {\mathcal{AP}}}^q (R)$}
\put(155,0){${\bf Gh}(R)$}
\put(215,55){\vector(-1,-1){43}}
\put(120,55){\vector(1,-1) {43}}
\put(165,67){$\theta_R^q$}
\put(165,54){$\cong$}
\put(205,35){$ \varphi_R^q$}
\put(120,34){$\tilde \varphi_R^q$}
\put(163,38){$\curvearrowright$}
\end{picture}
\end{equation}
is commutative.
The set of morphisms of
${\bf {\mathcal{AP}}}^q$ are same with that of ${\bf {\mathcal{AP}}}$
(see Section 4).

Summing up our arguments until now,
we obtain the following theorem:
\begin{thm}
The functors
$W^q$, $\mathcal C(F_q,\cdot)$, $\Delta^q$, and ${\bf {\mathcal{AP}}}^q$
are equivalent for all $q\in \mathbb Z.$
And each natural equivalence among them is compatible with
the Frobenius and Verschiebung operators.
\end{thm}
\small{

}
\end{document}